\documentclass[11pt]{amsart}
\usepackage{amssymb,amsmath,enumerate}
\numberwithin{equation}{section}
\theoremstyle{plain}
\newtheorem{thrm}{Theorem}[section]
\newtheorem{lemma}[thrm]{Lemma}
\newtheorem{prop}[thrm]{Proposition}

\newtheorem{rmrk}[thrm]{Remark}

\newtheorem{Question}[thrm]{Question}

\begin{document}
\begin{abstract}
Motivated by the equation satisfied by the extremals of certain
Hardy-Sobolev type inequalities, we show sharp $L^q$ regularity for
finite energy solutions of p-laplace equations involving critical
exponents and possible singularity on a sub-space of $\mathbb{R}^n$,
which imply asymptotic behavior of the solutions at infinity. In
addition, we find the best constant and extremals in the case of the
considered $L^2$ Hardy-Sobolev inequality.
\end{abstract}

\newcommand{\Om}{\Omega}
\newcommand{\Rm}{\mathbb{R}^n}
\newcommand{\Rnmk}{\mathbb{R}^{n-k}}
\newcommand{\Rmminusk}{\mathbb{R}^{n-k}}
\newcommand{\Rk}{\mathbb{R}^k}
\newcommand{\Rn}{\mathbb{R}^n}
\newcommand{\Hn}{\mathbb{H}^n}
\newcommand{\bG}{\mathbb{R}^n}

\newcommand{\dom}{{\mathcal{D}}\,^{1,p}(\Omega)}
\newcommand{\domG}{{\mathcal{D}}\,^{1,p}(\bG)}
\newcommand{\domO}{{\mathcal{D}}\,^{1,p}(\Omega)}
\newcommand{\domain}[1]{{\mathcal{D}}\,^{1,2}(#1)}
\newcommand{\pn}[1]{\norm{#1}_{\,\domG}}
\newcommand{\pnn}[2]{\norm{#1}_{\,L^{p}(#2)}}
\newcommand{\twost}{\frac {2}{n-2}}
\newcommand{\pst}{\frac {pn}{n-p}}
\newcommand{\psst}{\frac {p(n-s)}{n-p}}
\newcommand{\qn}[1]{\norm{#1}_{\,L^{p\text{*}}(G)}}
\newcommand{\qnn}[2]{\norm{#1}_{\,L^{p\text{*}}(#2)}}
\newcommand{\dl}[1]{\delta_{#1}}
\newcommand{\tr}[1]{\tau_{#1}}

\newcommand{\dis}[1]{\displaystyle #1}

\newcommand{\algg}{\mathfrak g}
\newcommand{\Lap}{\mathcal{L}}
\newcommand{\lap}{\triangle}
\newcommand{\gv}{\mathfrak v}
\newcommand{\gz}{\mathfrak z}
\newcommand{\norm}[1]{\lVert#1\rVert}
\newcommand{\abs}[1]{\lvert #1 \rvert}
\newcommand{\e}{\textbf {e}}

\newcommand{\st}{\text{*}}
\newcommand{\dxj}[1]{\frac{\partial {#1}}{\partial x_{j}}}
\newcommand{\dt}{\frac{d}{dt}}
\newcommand{\dxa}[1]{\frac {\partial {#1}} {\partial x_{\alpha}} }
\newcommand{\dr}[1]{\frac{\partial {#1}}{\partial r}}
\newcommand{\dyb}[1]{\frac {\partial {#1}} {\partial y_{\beta}} }
\newcommand{\dxap}[1]{\frac {\partial {#1}} {\partial x_{{\alpha'}}} }
\newcommand{\dybp}[1]{\frac {\partial {#1}} {\partial y_{{\beta'}}} }
\newcommand{\dxayb}[1]{\frac {\partial^2 {#1}} {\partial x_{\alpha}{\partial y_{\beta}}} }
\newcommand{\dxaxap}[1]{\frac {\partial^2 {#1}}{\partial x_{\alpha}{\partial y_{{\alpha'}}}}}
\newcommand{\dxaxa}[1]{\frac {\partial^2 {#1}}{\partial x_{\alpha}{\partial y_{\alpha}}}}
\newcommand{\dybybp}[1]{\frac {\partial^2 {#1}}{\partial y_{\beta}{\partial y_{{\beta'}}}}}
\newcommand{\R}{\mathbb{R}}
\newcommand{\Cn}{\mathbb{C}^n}
\newcommand{\dell}{\delta_{\lambda}}
\newcommand{\zeps}{\zeta_{\epsilon}}
\newcommand{\Ab}{\Bar{A}}
\newcommand{\nz}{\lvert z \rvert}
\newcommand{\nw}{\lvert w \rvert}
\newcommand{\zb}{\Bar{z}}
\newcommand{\wb}{\Bar{w}}
\newcommand{\bx}{\boldsymbol x}
\newcommand{\by}{\boldsymbol y}

\newcommand{\vx}{\boldsymbol x}
\newcommand{\vy}{\boldsymbol y}
\newcommand{\vz}{\boldsymbol z}
\newcommand{\tphi}{\tilde\phi}
\newcommand{\dx}[1]{\frac {\partial {#1}} {\partial x} }
\newcommand{\dxi}[1]{\frac {\partial {#1}} {\partial \xi} }
\newcommand{\dy}[1]{\frac {\partial {#1}} {\partial y} }
\newcommand{\deta}[1]{\frac {\partial {#1}} {\partial \eta} }

\newcommand{\dmu}[1]{\frac {dz} {\abs x^{#1}} }
\newcommand{\dmus}{\frac {dz} {\abs x^{s}} }

\newcommand{\ueps}{{u_\epsilon}}
\newcommand{\veps}{{v_\epsilon}}
\newcommand{\bueps}{{{\bar u}_\epsilon}}

\title[Hardy-Sobolev inequalities]{ $L^p$ estimates and asymptotic behavior for finite energy solutions of extremals
 to Hardy-Sobolev inequalities}
\date{June 6, 2006}
\author{Dimiter Vassilev}
\address[Dimiter Vassilev]{
University of California, Riverside\\
Department of Mathematics\\
Riverside, CA 92521} \email{dvassilev@math.ucr.edu}

\maketitle

\vskip.6truein
\section{Introduction}
\vskip.2truein

This paper has three goals - prove sharp $L^q$ regularity for
solutions of nonlinear p-laplacian equations involving critical
exponent and a singularity on a lower dimensional subspace,
establish sharp rate of decay in the case $p=2$, and finally
determine the extremals in a related $L^2$ Hardy-Sobolev inequality.
We indicate some other possible applications concerning stationary
cylindrical states of the Vlasov-Poisson system and non-completeness
of metrics with finite volume on some noncompact manifolds.

The organization of the paper is as follows. In Section 2 we shall
study the $L^p$ regularity and asymptotic behavior at infinity of
non-linear equations involving critical growth. In its simplest
form, our result is the following. For $1<p<n$, we let $p'=\frac
{p}{p-1}$ and $p^*=\frac {np}{n-p}$ be correspondingly the H\"older
conjugate and the Sobolev conjugate exponents. Suppose $V\in
L^{\frac {p^*}{p^* - p}}\,(\Rn)$, and if $p>2$  assume further
$V\geq 0$. Let $u$ be a weak nonnegative solution in $\Rn$, i.e. $u$
is of finite energy, cf. \eqref{weakVsol}, of the inequality
\begin{equation}
-\text{div}\, (|\nabla u|^{p-2} \nabla u)\ \leq \ Vu^{p-1}.
\end{equation}
Then, $u\in L^q (\Rn)$ for any $\frac {p^*}{p'}<q\leq\infty$, which
in general cannot be improved further, and for any $0<\theta<1$
there exists a positive constant $C_\theta$, such that,
\begin{equation}
u(z)\ \leq\ \frac {C_\theta}{1+|z|^{\theta \frac{n-p}{p-1}}},
\quad\quad z\in\Rn.
\end{equation}
 We shall actually consider
inequalities with a possible singularity on a subspace of $\Rn$ in
the right-hand side, which arise naturally in the study of the
extremals of certain Hardy-Sobolev inequalities, see
 further below and \eqref{e:maineq} for the exact setting.  Our results are contained in Theorems
\ref{boundedness}, \ref{T:belowp*}, \ref{T:Asymptotic1} and
\ref{T:asymptotic_beh p-lapl}.  The asymptotic behavior of solutions
to elliptic equations with critical non-linearity has been studied
extensively.  Generally speaking, the results and the analysis
depend on whether one assumes a priori that the solution has finite
energy or one is given a smooth solution,  see \cite{S1}, \cite
{S2}, \cite{SW}, \cite{E},\cite{GS} and \cite{Li}, to name a few. We
work under the assumption of finite energy. The paper extends the
results of Egnell \cite{E}, which correspond to Theorems
\ref{boundedness} and  \ref{T:Asymptotic1} without the possible
singularity on a submanifold of $\Rn$.

In Sections 3 we consider the case $p=2$. The main result here is
that, under some natural additional assumption, $u$ has the same
asymptotic as the fundamental solution of the laplacian at infinity,
see Theorem \ref{T:asymptotic_beh lapl}. The proof exploits the fact
that $\Rn$ has a positive Yamabe invariant and thus the method will
be applicable to cases when the ambient space does not have a
conformal transformation with the properties of the Kelvin
transform. Furthermore, we work in the setting in which $V$ could
have a singularity on a subspace, which renders results based on
radial symmetry inapplicable.

Section 4 contains some calculations which show how one can
construct explicit solutions of
 the Euler-Lagrange equation of certain Hardy-Sobolev inequalities.
 In the subsequent section we show that these solutions are
 extremals for which the best constant in the considered
 Hardy-Sobolev inequalities is achieved.

In Section 5, we determine the best constant in the Hardy-Sobolev
embedding Theorem involving the distance to a subspace of $\Rn$. In
order to explain the considered inequality we need a few notations,
which shall be used throughout the paper.  For $ p\geq 1$ we define
the space $ \domG$ as the closure of $C_o^{\infty}(\bG)$ with
respect to the norm
\begin{equation}\label{normDo1p}
\norm {u}_{\domG}\ =\ \left(\int_{\bG} |\nabla u|^p dz\right)^{1/p}.
\end{equation}

\noindent Let $n\geq 3$ and $2\leq k\leq n$. For a point $z$ in
$\Rm=\Rk\times\Rmminusk$ we shall write $z=(x,y)$, where $x\in\Rk$
and $y\in\Rmminusk$. The following Hardy - Sobolev inequality was
proven in Theorem 2.1 of \cite{BT}.

\medskip
\begin{thrm}[\cite{BT}]\label{BT}
Let $n\geq 3$,  $2 \leq k \leq n$, and $p, \ s$ be real numbers
satisfying $1 < p < n$, $0 \leq s \leq p$, and $s< k$. There exists
a positive constant $S_{p,s}=S(s,p,n,k)$ such that for all $u\in
D^{1,\, p}\,(\Rm)$ we have
\begin{equation}
\Bigl ( \int_{\Rm} \frac {|u|^{\frac {p(n-s)}{n-p}}}{|x|^s}\, dz\
\Bigr )^{\frac {n-p}{p(n-s)}} \leq\ S_{p,s}\, \Bigl ( \,\int_{\Rm}
|\nabla u|^p\, dz \,\Bigr )^{\frac {1}{p}}.
\end{equation}
\end{thrm}
\medskip

\noindent When $k=n$  the above inequality becomes the
Caffarelli-Kohn-Nirenberg inequality, see \cite{CKN}, for which  the
optimal constant $S_{p,s}$ was found in \cite{GY}. The case $p=2$
was considered earlier in \cite{O} and \cite{GMGT}, where the
inequality is written in an equivalent, but slightly different form.
If we introduce the numbers $\sigma\ =\ \frac 12\frac
{s(n-p)}{p(n-s)}$, hence $0\leq \sigma < 1$, and $p_\sigma\ =\ \frac
{n-p}{p(n-s)}$ the above inequality becomes
\begin{equation}\label{L2HS}
\Bigl ( \int_{\Rm} \frac {|u|^{p_\sigma}}{|x|^{\sigma p_\sigma}} \,
dz \Bigr )^{1/p_\sigma}\ \leq \ C \Bigl ( \int_{\Rm} |\nabla u|^2\,
dz \Bigr )^{1/2},
\end{equation}
where $C$ is a positive constant. In the case $p=2$ and $k=n$ the
sharp constant was computed in \cite{GMGT}. When $\sigma <1$ the
extremals were found by Lieb in \cite{L}, while when $\sigma =1$ we
have the classical Hardy inequality, which does not have extremal
functions.

The main result of Section 5 is the proof of the following Theorem.

\medskip
\begin{thrm}\label{T:bestconstantHardySobolev}
Suppose $n\geq 3$ and $2 \leq k \leq n$. There exists a positive
constant $K=K_{n,k,2}$ such that for all $u\in D^{1,\, 2}\,(\Rm)$ we
have
\begin{equation}\label{sharpHS}
\Bigl (\  \int_{\Rmminusk}\int_{\Rk} \frac {|u|^{\frac
{2(n-1)}{n-2}}}{|x|}\, dxdy\  \Bigr )^{\frac {n-2}{2(n-1)}}\ \leq\
K\, \Bigl ( \ \int_{\Rm} |\nabla u|^2\, dz \ \Bigr )^{\frac {1}{2}}.
\end{equation}
Furthermore, $K$ is given in \eqref{K} and the positive extremals
are the functions

\begin{equation}
v\ =\ \lambda^{-(n-2)} \,\Bigl ( \frac {4}{(n-2)^2}\Bigr )^{-\frac
{n-2}{2}} \, {K^{-(n-1)}}\, \Bigl (  \ (|x|+ \frac
{n-2}{4a\lambda^2} )^2 + |y-y_o|^2\ \Bigr )^{-\frac {n-2}{2}} ,
\end{equation}

\noindent where  $\lambda \ > \ 0$,  $y_o\in\Rmminusk$.
\end{thrm}

 \medskip
Let us note that the results of Section 2 can be applied to the
non-negative extremals of the general Hardy-Sobolev inequality of
Theorem \ref{BT}.

Finally, Section 6 contains some further simple applications, which
indicate the direction of some future investigations.

The author would like to thank Luca Capogna, Scott Pauls and Jeremy
Tyson for organizing the Workshop on Minimal surfaces, Sub-Elliptic
PDE's and Geometric Analysis, Dartmouth College, March, 2005.  The
current paper grew out of a conversation with Richard Beals and
Yilong Ni on an unrelated question of the author concerning certain
equations on Carnot groups, \cite{BY}, which lead to Section 4 of
this paper. Thanks are also due to Galia Dafni, Nicola Garofalo, Qi
Zhang and Congming Li for their interest in the work and valuable
conversations.

 \vskip.6truein
\section{Regularity and asymptotic of weak solutions}
\vskip.2truein

The goal of this section is to prove sharp $L^q$ regularity for
solutions of the considered equations, which would then lead to
bounds on the rate of decay at infinity.

We start with two definitions. Let  $p$ and $s$ be as in Theorem
\ref{BT} and denote by $ p\st(s) $ the Hardy-Sobolev conjugate
\begin{equation}\label{psst}
 p\st (s)\ = \psst
\end{equation}
 and by $p' $ the H\"{o}lder conjugate $\dis{
p'=\frac {p}{p-1}}$. There is another exponent, which will play an
important role. For any $s$ as above we define the exponent $r=r(s)$
to be the Holder conjugate of the exponent $r'=r'(s)$ defined by

\begin{equation}\label{rprime}
r'\ =\ \frac {p^*}{p^*(s)-p},\quad\text{thus}\quad r=\frac
{n}{n-p+s}.
\end{equation}
Notice that $1 \leq r$, $0\leq rs \leq p$, and furthermore we have
the identity

\begin{equation}\label{rsproperty}
rp\ =\ p^*(rs).
\end{equation}
However, in general $rs$ could be bigger than $k$, and thus due to
the restriction $s<k$ in the Hardy-Sobolev inequality we shall
consider only the case
\begin{equation}\label{s-bound}
s(n-k)\ <\ k(n-p),
\end{equation}
which implies $rs<k$.
\medskip

\begin{thrm}\label{boundedness}
Let $\Om$ be an open subset of $\Rn$, which is not necessarily
bounded, $1< p <n$, $0 \leq s \leq p$, $s< k$ and $s(n-k)\ <\
k(n-p)$. Let $u\in \dom$ be a non-negative weak solution of the
inequality
\begin{equation}\label{e:maineq}
-\ \text{div}\, (|{\nabla}u|^{p-2}{\nabla}u)\ \leq\ V\frac {\
|u|^{p-2}}{|x|^s}u \quad \quad \text{in} \quad \Omega,
\end{equation}
i.e.,
\begin{equation}\label{weakVsol}
\int_{\Omega} \lvert {\nabla}u \rvert \,^{p-2} <{\nabla}u,
{\nabla}\phi > \,dz\ \leq \ \int_{\Omega} V\frac {\
|u|^{p-2}}{|x|^s}u\,\phi\,dz ,
\end{equation}
for every $0\leq \phi \in C^{\infty}_o (\Omega)$.

a) If $V\in L^{r'}(\Omega)$, then $u\in L^q\, (\dmu t)$ for any
$0\leq t<\min \{p,s\}$ and $q\geq p^*(s)$. In particular $u\in
L^q(\Om)$ for every $p^*\ \leq\ q\ <\ \infty.$

b) If $V\in L^{t_o}(\Om)\cap L^{r'}(\Omega)$ for some  $t_o\ >\ r'$,
then $u\in L^\infty(\Om)$.

\end{thrm}

\begin{rmrk}\label{R:L-infinity for sols}
a) As usual, if we have equality instead of inequality in the above
Theorem, the conclusion holds for any weak-solution, without a sign
condition. The proof below works with very minor changes. This is
the reason we use $|u|$ rather than  $u$ when we are dealing with a
non-negative function.

b) The use of the weighted $L^q$ spaces in part (a) of the Theorem
is essential.
\end{rmrk}
\medskip
\medskip

\begin{proof}

The assumption $V\in L^{r'}(\Om)$ together with the Hardy-Sobolev
inequality shows that \eqref{weakVsol} holds true for any
$\phi\in\dom$. This can be seen by approximating in the space $\dom$
by a sequence of test functions $\phi_n \in C^{\infty}_o (\Om)$,
which will allow to put the limit function in the left-hand side of
\eqref{weakVsol} as $\abs{\nabla u}^{p-1}\in L^{p'}$.  On the other
hand, for $\phi\in C^{\infty}_o (\Om)$, using the H\"older and
Hardy-Sobolev inequalities we have the estimate

\begin{align}
\int_\Om  |V|    \frac { |u|^{p-1}} {|x|^s}   \phi \, dz\ & \leq\
\Bigl( \int_\Om |V|^{r'} \Bigr)^{1/{r'}}\,
 \Bigl( \int_\Om \frac {|u|^{r(p-1)}} {|x|^{rs}} \phi^r\, \Bigr)^{1/r}
\\
 & \leq\  \Bigl ( \int_\Om |V|^{r'} \Bigr )^{1/{r'}} \, \Bigl [ \Bigl ( \int_\Om \frac {
|u|^{rp'(p-1)}}{|x|^{rs}}\, \Bigr )^{\frac {1}{rp'}}\, \Bigl (
\int_\Om \frac {\phi^{pr}}{|x|^{rs}}\,
\Bigr )^{1/p}\Bigr ]^{1/r}\\
& \leq\ S_p\ \norm {V}_{L^{r'}}\, \Bigl ( \int_\Om\frac {
|u|^{rp'(p-1)}}{|x|^{rs}}\,  \Bigr )^{\frac {1}{rp'}}\, \norm
{\nabla\phi}_{L^p} \quad \text{( using \eqref{rsproperty} )}\\
& \leq\ S_pS^{p-1}_{p,rs}\norm {V}_{L^{r'}}\, \norm{\nabla
u}^{p-1}_{L^p}\, \norm {\nabla\phi}_{L^p},
\end{align}
which allows to pass to the limit in the right-hand side of
\eqref{weakVsol}.

We turn to the proofs of a) and b).

a) Let $G(x)$ be a piece-wise smooth, globally Lipschitz function,
on the real line, and set

\begin{equation}\label{Fdef}
F(u)\ =\ \int^{u}_0 \abs {G'(t)}^p\, dt.
\end{equation}

\noindent Clearly, $F$ is a non-negative differentiable function
with bounded and continuous derivative. From the chain rule, $G(u),\
F(u)\in\dom$. In particular, $F(u)$ is a legitimate test function in
\eqref{weakVsol}. We are going to show that if $q$ is a number
$q\geq p^*(s)$ and $0\leq t<\min \{p,s\}$, then $u\in L^q\, (\dmus)$
implies $u\in \dom\cap L^{\kappa_t q}\, (\dmu t)$, where
$\kappa_t=\frac {p^*(t)}{p}$, and for some positive constant $C$
depending on $q$ we have
\begin{equation}\label{Lqclaim}
 \norm u _{L^{\kappa_t q}\, (\dmu t)} \
\leq\ C \ \norm u_{L^q\, (\dmus)}.
\end{equation}
Notice that we can apply the Hardy-Sobolev inequality replacing the
exponent $s$ with the exponent $t$. Furthermore, we require $t<p$ as
then we have

\begin{equation}\label{kappat}
\kappa_t\ =\ \frac {p^*(t)}{p} \ >\ 1.
\end{equation}

With $\phi=F(u)$,  taking into account $F'(u)\ =\ \abs{G'(u)}^p$,
the left-hand side of \eqref{e:maineq} can be rewritten as

\begin{equation}\label{LHS}
\int_{\Om} \lvert {\nabla}u \rvert \,^{p-2} <{\nabla}u, {\nabla}
F(u) > \,dH\ =\ \int_{\Om} \abs{{\nabla} G(u)}^p.
\end{equation}

\noindent For $q\geq p^*(s)$, hence $q\geq p^*(p)=p$, we define the
function $G(t)$ on the real line in the following way,

\begin{equation}\label{e:definition of G}
G(t)=\begin{cases}
     sign\,(t)\,\abs{t}^{\frac qp} & \quad\quad \text{ if } 0\ \leq\ \abs{t}\ \leq l, \\
     l^{\frac qp -1}t & \quad\quad \text{ if } l\ <\ \abs{t}.
     \end{cases}
\end{equation}

\noindent From the power growth of $G$, besides the above
properties, this function satisfies also

\begin{equation}\label{FGineq}
\abs{u}^{p-1}\abs{F(u)}\ \leq\ C(q)\abs{ G(u)}^p \ \leq\ C(q)
\abs{u}^q.
\end{equation}

\noindent The constant $C(q)$ depends also on $p$, but this is a
fixed quantity for us. At this moment the value of $C(q)$ is not
important, but an easy calculation shows that $C(q)\ \leq\ C\,
q^{p-1}$ with $C$ depending on $p$. We will use this in part b). For
$ M>0 $ to be fixed in a moment we estimate the integral in the
right-hand side of \eqref{weakVsol} as follows.

\begin{multline}\label{RHS}
\int_{\Om} |V|\frac {\abs{u}^{p-1}}{\abs{x}^s} F(u)\, dz  \\
\\
=\  \int_{(\abs{V} \leq M)} \abs V\frac {\abs{u}^{p-1}}{\abs{x}^s} F(u)\, dz \  +\ \int_{(\abs{V}>M)} \abs V\frac {\abs{u}^{p-1}}{\abs{x}^s} F(u)\, dH \\
\\
                               \leq \ C(q) \int_{(\abs{V} \leq M)} \abs {V}\frac {\abs{G(u)}^p}{\abs {x}^s}\, dz\
                               + \ C(q){\Bigl( \int_{(\abs{V} >M)} \abs{V}^{r'}
                               \Bigr) }^ {\frac {1}{r'}}
{ \Bigl(  \int_{\Om}\frac {\abs{G(u)}^{pr}}{\abs {x}^{sr}} \, dz
\Bigr)}^{\frac {1}{r}}\\
\\
\leq \ C(q) \int_{(\abs{V} \leq M)} \abs {V}\frac {\abs{u}^q}{\abs
{x}^s}\, dz\
                               + \ C(q)S_{p,rs}^p{\Bigl( \int_{(\abs{V} >M)} \abs{V}^{r'}
                               \Bigr) }^ {\frac {1}{r'}}
{ \Bigl(  \int_{\Om}\abs{\nabla G(u)}^{p} \, dz \Bigr)}\\
\\
\leq \ C(q) M \ \norm {u}^q_{L^q\, (\dmus)}\
                               + \ C(q)S_{p,rs}^p \norm {V}_{L^{r'}\,
                               (|V|>M)}\
\norm {\nabla G(u)}^p_{L^p}.
\end{multline}
At this point we fix once and for all the constant $M$, so that

\begin{equation*}
C(q)S_{p,rs}^p{\Bigl( \int_{(\abs{V} >M)} \abs{V}^{r'}\, dH \Bigr)
}^ {\frac {1}{r'}}\ \leq\ \frac{1}{2},
\end {equation*}

\noindent which can be done because $ V \in L^{r'} $. Putting
together \eqref{LHS} and \eqref{RHS}, and using the Hardy-Sobolev
inequality, we come to

\begin{equation}
\norm {G(u)}^p_{L^{p^*(t)}(\dmu t)}\ \leq \ C(q) M \ \norm
{u}^q_{L^q\, (\dmus)}
\end{equation}

\noindent By Fatou's theorem we can let $l$ in the definition of $G$
to infinity and obtain

\begin{equation*}
\norm {u}^q_{L^{\kappa_t q}(\dmu t)}\ \leq \ C(q) M \ \norm
{u}^q_{L^q\, (\dmus)}
\end{equation*}
The proof of a) is finished.

b) Let us observe that the assumption $t_o \
> \ r'$ implies that $t_o'\ < \ r$ and thus $0\ <\ t_o' s\ <\
rs$. Therefore, for any $q\geq p^*(s)$ the norm $\norm{u}_{L^{q
t_o'}(\dmu {st_o'})} $ is finite from part (a). We are going to
prove that the $L^q (\dmu {t_o' s})$ norms of $u$ are uniformly
bounded by the $L^{q_o}(\dmu {t_o' s})$ norm of $u$,  where $q_o\ =\
t_o'\,p^*(s)$. We shall do this by iteration and find a sequence
$q_k$ which approaches infinity as $k\rightarrow\infty$.

Let $q\geq p^*(s)$. We use again the function $F(u)$ from part (a)
in the weak form \eqref{weakVsol} of our equation. The left-hand
side is estimated from below as before, see \eqref{LHS}. This time,
though, we use H\"{o}lder's inequality to estimate from above the
right-hand side,

\begin{multline}\label{RHSt}
\int_{\Om} \abs V\abs{u}^{p-2}u F(u)\, dH\ \leq\ \norm{V}_{L^{t_o}}\, \norm{\abs{u}^{p-1}F(u)}_{L^{t_o'}}\\
\leq\ \norm{V}_{t_o}\, \norm{C(q)\frac
{\abs{G(u)}^p}{|x|^s}}_{t_o'}\ \leq\ C(q) \norm{V}_{L^{t_o}}\,
\norm{u}_{L^{q t_o'}(\dmu {st_o'})}^q.
\end{multline}

\noindent With the estimate from below we come to

\[
\norm{\nabla G(u)}_{L^p} ^p\ \leq\ C(q)\,\norm{V}_{L^{t_o}}\,
\norm{u}_{L^{q t_o'}(\dmu {st_o'})}^q.
\]

\noindent Using the Hardy-Sobolev inequality and then letting $l\,
\rightarrow\, \infty$ we obtain

\begin{equation}\label{iterate}
\norm {u}^q_{L^{q\frac {p^*(st_o')}{p}}(\dmu {st_o'})}\ \leq\
C\,C(q)\,\norm{V}_{L^{t_o}}\, \norm{u}_{L^{q t_o'}(\dmu {st_o'})}^q,
\end{equation}
where $C$ is independent of $q$.

 \noindent Let $\dis{\delta\ =\
\frac {p^*(st_o')}{p\,t_o'}}$. A small calculation shows that
$\delta > 1$ exactly when $t_o>r'$. With this notation we can
rewrite \eqref{iterate} as

\begin{equation}\label{iter}
\norm {u} _{L^{\delta q t_o'}(\dmu {st_o'})}^q\ \leq\
\Bigl[C\,C(q)\Bigr]^{\frac 1q}\,\norm{V}_{L^{t_o}}^{\frac 1q}\,
\norm{u}_{L^{q t_o'}(\dmu {st_o'})}^q.
\end{equation}

\noindent Recall that $C(q)\ \leq\ C q^{p-1}$. At this point we
define $q_o\ =\ p^*(s) t_o'$ and $q_k\ =\ \delta^k q_o$, and after a
simple induction we obtain

\begin{equation}
\norm{u}_{L^{q_k}(\dmu {st_o'})}\ \leq\  \Bigl \{ \prod^{k-1}_{j=0}
\bigl[C\, q_{j}^{p-1}\bigr]^{\frac {1}{q_j}}\ \Bigr \}
\norm{V}_{L^{t_o}} ^{\sum^{k-1}_{j=0} \frac {1}{q_{j}}}\
\norm{u}_{L^{q_o}(\dmu {st_o'})}.
\end{equation}

\noindent Let us observe that the right-hand side is finite,

\begin{equation}\label{e:finitesums}
\sum^{\infty}_{j=0} \frac {1}{q_{j}}\ =\ \frac
{1}{q_o}\sum^{\infty}_{j=1} \frac {1}{\delta^{j}} <\infty
\quad\quad\text{ and }\quad\quad \sum^{\infty}_{j=1} \frac {\log
q_{j}}{q_{j}} \ <\ \infty,
\end{equation}

thanks to $\delta\ >\ 1.$ Letting $k \rightarrow \infty $ we obtain

\[
\norm{u}_{\infty}\ \leq\ C\, \norm{u}_{L^{q_o}(\dmu {st_o'})}.
\]

\end{proof}

\medskip
\begin{rmrk}\label{r:localboundedness}
We should keep in mind that the local version of Theorem
\ref{boundedness} is also valid. In other words,
 we can replace all the spaces in the statement of the Theorem with
 their local version. The proof is accomplished in a very similar
 fashion by introducing a local cut-off function.
\end{rmrk}

\medskip

With the above Theorem we turn to the equation, in fact slightly
more general equation, satisfied by the extremals of the
Hardy-Sobolev inequality.

\medskip

\begin{thrm}\label{T:L-infinity}
Let $\Om$ be an open subset of $\Rn$, which is not necessarily
bounded, $1< p <n$, $0 \leq s \leq p$, $s< k$ and $s(n-k)\ <\
k(n-p)$. If $R\in L^\infty$ and $u \in\dom$ is a weak non-negative
solution to equation

\[
-\ \text{div}\, (|{\nabla}u|^{p-2}{\nabla}u)\ \leq \  R(z)\frac {\
|u|^{p^*(s)-2}u}{|x|^s}u \quad \quad \text{in} \quad \Omega,
\]
then $u\in L^\infty(\Om)$.
\end{thrm}

\medskip

\begin{proof}
We define $V\ =\ R\abs{u}^{p^*(s) - p}$. From the Hardy-Sobolev
inequality we have $u\in L^{p^*(s)} (\Om)$ and thus $V\in L^{\frac
{p^*(s)}{p^*(s) - p}} (\Om)$. Since $\dis {r'\ =\ \frac
{p^*(s)}{p^*(s) - p}}$, part (a) of Theorem \ref{boundedness} shows
that $u\in L^q\,(\Om)$ for $p^*\leq q < \infty$. Therefore $V\in
L^{\frac {q}{p^* - p}} (\Om)$ for any such $q$ and thus by part (b)
of the same Theorem we conclude $u\in L^\infty(\Om)$.
\end{proof}

\medskip

In the next Theorem we show that one can lower the exponent $p^*$ in
the $L^q$ regularity of $u$, when $s=0$. A similar result in the
case $p=2, \ s=0$ and $R=\abs{u}^{p^*-p}$ was achieved in \cite{LU},
see also \cite{BK} and \cite{GL}. We have to overcome some
complications due to the more general structure of the equation, the
possible singularity, and the lack of monotonicity of the considered
operator.

\medskip
\begin{thrm}\label{T:belowp*}
Let $\Om$ be an open subset of $\Rn$, which is not necessarily
bounded, $1< p <n$, $0 \leq s \leq p$, $s< k$ and $s(n-k)\ <\
k(n-p)$. Suppose $R\in L^{r'}$ and $V_o\in L^1\cap L^{r'}$, and in
the case $p>2$ assume $R$ and $V_o$ are non-negative, $R,\ V_o\geq
0$. If $u$ is a non-negative locally bounded weak solution of the
equation

\begin{equation}\label{e:main_ineq+Vo}
-\triangle_p\,  u \ \leq \ R\frac {\abs{u}^{p-2}}{\abs {x}^s}u\ +\
V_o,
\end{equation}
 then $u\in L^q$ for every $ \frac
{p^*}{p'}\ < q\ \leq \ p^*$.
\end{thrm}

\begin{rmrk} The condition that $u$ is locally bounded holds for
example when $V_o=0$ and $R\in L^{r'}\cap L^{t_o}$, for some $t_o
>r'$ by Theorem \ref{boundedness}.
\end{rmrk}
\medskip

\begin{proof}
Let $0<\theta <1/p$ be arbitrarily fixed. Our task then is to show
that the function $u^{1-\theta}\in L^{p^*}$.

In the first part of the proof we shall exploit the variational
structure of the problem in order to construct suitable test
functions, which shall be used in the second part of the proof.
Suppose $V, g \ \in L^1\cap L^{r'}$ are two given functions.

Consider the functional
\begin{equation}
E(v)\ =\ \frac 1p\, \int_\Om \abs {\nabla v}^p\, dz\ -\ \frac 1p\,
\int_\Om V \frac {\abs {v}^p}{\abs{x}^s}\, dz\ -\ \int_\Om g v\, dz.
\end{equation}
Note that $E$ is coercive when $\norm {V}_ {   L^{r'} (\Om)  }  $ is
small. Indeed, we have
\[
E(v) \ \geq\ \frac 1p \, \norm {v}^p _{\domO} \ -\  S_{p, rs}
^{p-1}\
 \frac 1p\, \norm {V}  _{  L^{r'} }\ \norm {v}^p_{\domO} \ -\ S_{p,
0}\ \norm {g} _{  L ^{(p^*)'} }\ \norm {v}_{\domO},
\]
taking into account \eqref{e:p*'vs r'} to justify the finiteness of
the norm of $g$. Furthermore, $E$ is weakly lower semi-continuous
and $\domO$ is a reflexive Banach space. Therefore, provided that
$V$ has a small norm, there exists a minimizer, which is a solution
of
\begin{equation}\label{e:u bar eqn}
-\triangle_p \ v \ = \  V \frac {\abs{v}^{p-2}}{\abs{x}^s}\, v\ +\
g.
\end{equation}
Moreover, any solution of the above equation satisfies
\begin{equation}\label{e:uniform bound for uk}
\norm {v}_{\dom}\ \leq\  S_p^{1/(p-1)}\
(1-S_{p,rs}^p\norm{V}_{L^{r'}} )^{1/(p-1)}\
\norm{g}^{1/(p-1)}_{L^{(p^*)'}}.
\end{equation}

 With this in mind, suppose $\epsilon$ is a
given positive constant. Since  $R\in L^{r'}$ we can fix a large
$R_o>0$ and a small $\delta>0$, such that
\begin{equation}\label{e:choose Ro}
\int_{ \Rn\setminus B_{R_o} } \abs{R}^{\, r'}\, dz\ \leq\ \frac 12\
\epsilon \quad \text{ and }\quad
 \int_{ \{ \abs{x}<2\delta\} } \abs{R}^{\, r'}\, dz\
\leq\ \frac 12\ \epsilon.
\end{equation}
Let $\alpha\in C^\infty$ be a function, $0\leq \alpha\leq 1$, with
$$
\text{supp}\ \, \alpha \subseteq \{ {\abs{x}<2\delta}\} \cup \{
\Rn\setminus B_{R_o} \},
 \hskip.3truein\alpha\equiv 1\quad\text{on}\quad \,
\{ {\abs{x}<\delta} \} \cup \{ \Rn\setminus B_{2R_o} \}.
$$

\noindent In particular supp $ (1-\alpha)\subset B_{2R_o}\cap \{
{\abs{x}>\delta} \} $, and hence due to the local boundedness of $u$
we have $$g=V_o \ +\ (1-\alpha)\ R(z)\ \frac {\abs{u}^{p -
2}}{\abs{x}^s}\ u\in\ L^{1}\cap L^{r'}.$$

\noindent For every $k\in\mathbb{N}$ let $\alpha_k\in C^\infty_o\,
(\Rn)$ be a function, $0\leq \alpha_k\leq 1$,  satisfying
\[
\text{supp}\ \alpha_k\ \subset B_{2^{k+1}R_o}\setminus
 \{ \abs{x}<\frac {\delta}{2^{k+1}} \} \quad \text{ with } \quad \alpha_k=1\quad \text{  on } \quad
B_{2^{k}R_o}\setminus
 \{ \abs{x}<\frac {\delta}{2^{k}} \}.
\]
Notice that $\alpha_k \nearrow 1$ a.e. as $k\rightarrow\infty$.
Define $V = \alpha\, R(z)$ and  $V_k = \alpha\, \alpha_k\, R(z)$.
Using the properties of the cut-offs  we see that these functions
enjoy the following properties
\[
V\in L^{r'},\quad\quad V_k \ \in\  L^{1}\cap L^{r'}\quad\text{ as
}\quad \text{supp}\ V_k\subset\subset\Rn, 
\]
and
\[
V_k \nearrow V \ \text{as}\ k\rightarrow\infty.
\]
In addition, since
\[
\int_\Om \abs {V}^{r'}\, dz\ \leq \  \int_\Om \alpha\ \abs {R}^{\,
r'} \, dz\ \leq \  \epsilon,
\]
the functions $V$ and hence $V_k$ have small $L^{r'}$ norms, which
can be made less than $\epsilon$ by taking $R_o$ sufficiently large
and $\delta$ sufficiently small  in \eqref{e:choose Ro}. From now on
we assume that $R_o$ and $\delta$ have been fixed in the above
described manner so that $E$ is coercive and further so that $\norm
{V}_{L^{r'}}\ \leq \ \epsilon$ with $\epsilon$ to be fixed later
independently of $k$, and in fact depending only on $p,\ u, \ \norm
{R}_{L^\infty}$ and the Sobolev constant.

Rewriting the equation given in the Theorem, we have that $u$ is a
given non-negative solution of
\begin{equation}
-\text{div}\, (|{\nabla}u|^{p-2}{\nabla}u)\ \leq \ \  V  \frac {\
|u|^{p-2}}{\abs{x}^s}\, u\ +\ g.
\end{equation}
For every $k\in\mathbb{N}$ let $u_k$ be a solution of
\begin{equation}\label{e:equation for uk}
-\text{div}\, (|{\nabla}u_k|^{p-2}{\nabla}u_k)\ = \ \  V_k \frac {\
|u_k|^{p-2}}{\abs{x}^s}\,u_k\ +\ g.
\end{equation}

\noindent Moreover, when $p>2$ we assume that $R,\ V_o\geq 0$ and
thus if we take $u_k$ to be a minimizer of
\begin{equation}
E_k(v)\ =\ \frac 1p\, \int_\Om \abs {\nabla v}^p\, dz\ -\ \frac 1p\,
\int_\Om V_k \frac {\abs {v}^p}{\abs{x}^s}\, dz\ -\ \int_\Om g v\,
dz.
\end{equation}
then $u_k$ is  non-negative, i.e., when $p>2$ we have $u_k \geq 0$.
This shall be used at the very end of the proof.

 Next we define the necessary cut-off, which shall be used
in the final step. Let $\eta_m(t)$ be the following function
\begin{equation}
\eta_m (t)\ = \
    \begin{cases}
    t, &   t>1/m   \\
    & \\
      m^{\frac {\theta p}{1-p\, \theta} }  \ t^{\frac {1}{1-p\, \theta} }, &  0\leq\ t \leq \frac
    {1}{m}.
    ,
    \end{cases}
\end{equation}
Observe that $\eta_m$ is a continuous function and
\[
0\leq \eta'_m \leq \max \{ 1, \frac {1}{1-p\, \theta} \}\ =\ \frac
{1}{1-p\, \theta},
\]
which implies the useful fact
\[
\eta'_m(t)\ \leq\ \frac {1}{1-p\, \theta}\ t.
\]

\noindent We define also
\begin{equation}
\phi_m(t)\ =\ \eta_m^{1-p\theta}\ \text{ and } \ f_m(t)\ =\
\eta_m^{1-\theta}.
\end{equation}
A short calculation gives
\begin{equation}
\phi_m (t)\ = \
    \begin{cases}
    t^{1-p\theta}, &   t>1/m   \\
    & \\
      m^{p\theta }t, &   t \leq \frac
    {1}{m}
    \end{cases}
\quad\quad  \phi'_m (t)\ \leq \
    \begin{cases}
    {(1-p\, \theta)}\ m^{p\theta}, &   t>1/m   \\
    & \\
      m^{p\theta }, &   t \leq \frac
    {1}{m}
    \end{cases}
\end{equation}

\noindent In particular, for every fixed $m$,  we have that
$\phi'_m$ is a bounded function and thus if  $0\leq\ v\in\dom$ then
$0\leq \phi_m (v)\in\dom $. Since $1-\theta\
>\ 1-p\theta$ the derivative $f'_m(t)$ is also bounded and thus $f_m(v)\,
\in\dom$.  From now on, for simplicity, given a function $0\leq\,
v\in\dom$ we let

$$\eta=\eta_m(v), \hskip.5truein
\phi=\phi_m(v)\hskip.25truein\text{ and }\hskip.25truein f=f_m(v),
$$
all of which, due to the chain rule,  are functions from
 $\dom$.

A small but very important calculation shows that we have
\begin{equation*}
\abs {\nabla f}^p\ = \ (1-p\,\theta)^{p-1}\Bigl ( \frac
{1-\theta}{1-p\theta} \Bigr )^p\, \abs {\nabla \eta}^{p-2}\, \nabla
\eta\cdot \nabla \phi .
\end{equation*}

\noindent With $c_\theta\ =\ \Bigl ( \frac {1-\theta}{1-p\,\theta}
\Bigr )^p $, using  $\nabla v\cdot\nabla\phi \geq 0$ and the above
identity, we compute

\begin{multline*}
\int_{\Om} \abs{\nabla \eta}^{p-2}\  \nabla \eta\cdot \nabla \phi \,
dz\ =\ \int_{\{v<\frac 1m\}} \abs{\nabla \eta}^{p-2} \ \nabla \eta
\cdot \nabla \phi \, dz\ + \ \int_{ \{ \frac {1}{m}<v \} }
\abs{\nabla \eta}^{p-2}\ \nabla \eta\cdot \nabla \phi \, dz  \\
\\
= \ \frac     {   m^{  (p-1) p\, \theta/(1 - p\, \theta) }   }
             {           (1-p\,\theta)^{p-1}              }
 \    \int_{\{v<\frac 1m\}} v^{ (p-1)p\, \theta /(1-p\, \theta)}\
\abs{\nabla u}^{p-2} \nabla v \cdot \nabla \phi \,  dz\\
\\
 + \ \int_{\{\frac {1}{m}<v\}}
\abs{\nabla v}^{p-2} \  \nabla v \cdot \nabla \phi \, dz    \\
\\
\leq\ \frac {1}{(1-p\, \theta)^{p-1}}  \  \int_{\{v<\frac 1m\}}
  \abs{\nabla u}^{p-2}  \   \nabla v
\cdot \nabla \phi \, dz\ + \ \int_{\{\frac {1}{m}<v\}}
\abs{\nabla v}^{p-2}  \  \nabla v\cdot \nabla \phi \, dz    \\
\\
\leq \int_\Om \abs {\nabla v} ^{p-2} \ \nabla v\cdot \nabla \phi \,
dz  \quad \text{as}\quad  1\ \leq \ 1/(1-p\, \theta)
\end{multline*}
Therefore the following bound holds
\begin{equation}
 \int_{\Om} \abs {\nabla f}^p \, dz\ \leq\ c_\theta  \int_\Om \abs {\nabla v} ^{p-2} \ \nabla v\cdot \nabla \phi \,
dz .
\end{equation}

\noindent Let us set  \[ v=u_k^+ \] and use $\phi$ as a test
function in \eqref{e:equation for uk}. The inequality above gives an
estimate of the left-hand side. For the first term on the other side
we have

\begin{align*}
\int_\Om V_k \frac {\abs{v}^{p-1}}{\abs{x}^s}\ {\phi}\, dz\ & =\
\int_{\{ v>1/m \}} V_k \frac
 {\abs{v}^{p-1}}{\abs{x}^s}\ {\phi}\, dz \
 +\ \int_{\{v< \frac {1}{m}\}} V_k \frac {\abs{v}^{p-1}}{\abs{x}^s}\  {\phi}\, dz \\
&\\
 &\leq\ S_p^p\ \int_{\Om} V_k \frac
{\abs{v}^{p-1+1-\theta p}}{\abs{x}^s}\, dz\ +\ m^{\theta
p}\int_{\{v< \frac
{1}{m}\}} V_k \frac {\abs{v}^{p}}{\abs{x}^s}\, dz \\
&\\
& \leq\ S_p^p \norm {V_k}_{  L^{r'} }\int_{\Om} \abs {\nabla
f}^p \, dz\ +\  m^{\theta p - p}\ \frac {2^{s(k+1)}}{\delta^s}\ \norm {V_k}_{L^{1}} \\
&\\
& \leq\  S_p^p \norm {V}_{  L^{r'} }\int_{\Om} \abs {\nabla f}^p \,
dz\ +\ m^{\theta p - p}\ \frac {2^{s(k+1)}}{\delta^s}\ \norm
{V_k}_{L^{1}}
\end{align*}
using $pr=p^*(rs)$ and the fact that $f\in\dom$.  The other term on
the right-hand side can be estimated in the following way

\begin{multline}
\int_\Om g \phi\, dz\ \leq \ \int_{\{v<\frac 1m\}} \abs {g}\ \phi
\,dz\ +\ \int_{\{\frac {1}{m}<v\}} \abs {g}\ \phi \, dz\\
\\
\leq\ m^{p\,\theta}\int_{\{v<\frac 1m\}} \abs {g}\  v\,dz\ +\
\int_{\{\frac {1}{m}<v\}} \abs {g} \ v^{1-\theta
p} \, dz\\
  \\
  \leq\ m^{\theta
p - 1} \norm{g}_{L^1}\ +\
 \norm{g}_{L^q}\ \norm{v}^{p^*/q}_{L^{p^*}},
\end{multline}

\noindent where
\begin{equation}\label{e:definition of t}
q \ =\ \Bigl ( \frac {p^*}{1-p\theta}\Bigr )'\ =\ \frac
{p^*}{p^*-(1-p\theta)}\ >\ 1,
\end{equation}
as $1-p\, \theta>0$. Recall we are assuming only $g \in L^1\cap
L^{r'}$. The definition \eqref{e:definition of t} of $q$ shows that
$q < (p^*)'$, while
\begin{equation}
p^*(p)=p\ \leq\ p^*(s)\ \leq\ p^*(0)=p^*
\end{equation}
and thus
\begin{equation}\label{e:p*'vs r'}
r'=\ \frac {p^*}{p^*(s)-p}\ \geq\ \frac{p^*}{p^*-p}\ =\ \bigl (p^*/p
\bigr )' \
> \  \frac {p^*}{p^*-1}\ =\ (p^*)'.
\end{equation}
Therefore we have
\[
q\ <\ r'
\]
and $g\in L^q_\text {loc}$ for every $1\leq q < r'$.

Putting the above three estimates together we have shown that the
gradient of $f$ satisfies the inequality ( with $S_p\ = \ S_{p,\,
o}$ )
\begin{multline}
\frac {1}{c_\theta} \int_{\Om} \abs {\nabla f}^p \, dz\ \leq \
S^p_{p} \norm {V}_{L^{r'}} \int_{\Om} \abs{\nabla f}^{p} \, dz\ +\ m^{\theta p - p}\ \frac {2^{s(k+1)}}{\delta^s}\ \norm {V_k}_{L^{1}}\\
\hskip3.5truein\\
+\ m^{\theta p - 1} \norm{g}_{L^1} +\
 \norm{g}_{L^q(v>1/m)}\ \norm{v}^{p^*/q}_{L^{p^*}(v>1/m)}.\quad\quad
\end{multline}
\medskip

\noindent With the above estimates at hand, we can conclude by
moving the first term from the right side to the left, and then
letting $m\rightarrow \infty$, provided we have

\begin{equation}\label{e:needed small norm of V}
    1-c_\theta S^p_{p}\norm {V}_{L^{r'}}>0, \text{ i.e.}, \quad \norm {V}_{L^{r'}}\ll
\end{equation}

\noindent Taking the initial $R_o$ sufficiently large at the start
of the proof, we can let $m\rightarrow\infty$, recall
$f_m\rightarrow (v^+)^{1-\theta}$, which gives the inequality
\begin{equation*}
\int_{\Om} \abs {\nabla (v^+)^{1-\theta}}^p \, dz\ \leq \
C_{p,\theta}
 \norm{g}_{L^q(\Om)}\ \norm{v}^{p^*/q}_{L^{p^*}(\Om)},
\end{equation*}

Working similarly with $v^{-}$ we can prove eventually, recall
$v=u_k$, that the $\domG$ norms of $u^{1-\theta}_k$ satisfy
\begin{equation}
\norm {u^{1-\theta}_k}^p_{\domG}\ =\ C_{p,\theta}\,
 \norm{g}_{L^q(\Om)}\ \norm{u_k}^{ {p^*/q}}_{L^{p^*}(\Om)}
\end{equation}
and hence they are uniformly bounded in view of \eqref{e:uniform
bound for uk}.

The proof is finished by letting $k\rightarrow\infty$, but the
argument in the cases $p=2$ and  $p>2$ are different. When $p>2$ we
take into account that by construction $u_k\geq 0$, cf. the line
after \eqref{e:equation for uk}, and if $u_k\rightarrow \bar u$
weakly in $\dom$, where
\[
-\triangle_p \ \bar u \ = \  V \frac {{\bar u}^{p-1}}{\abs{x}^s} \
+\ g,
\]
from Lemma \ref{L:uniqueness p-lapl} we conclude $u^{1-\theta}\in
L^{p^*}$. If $p=2$ the situation is simpler since we can use the
monotonicity. In fact, since $\norm {V}_{L^{r'}}$ is small and $u, \
\bar u\in \, \dom$ H\"older's inequality and the strong monotonicity
of the laplacian
 give $(u - \bar u)^+ = 0$.
Indeed, using $w=(u - \bar u)^+ $ as a test function in the
inequality
\[
-\triangle u \ +\ \triangle \bar u\ \leq \ \ V \frac {u- \bar
u}{\abs{x}^s}
\]
we see that, cf. \eqref{rsproperty},
\begin{equation*}
\norm {w}^2_{\domain \Om} \leq C\ \norm{V}_{L^{r'}}\  \Bigl ( \int
\frac {w^{2r}}{ \abs{x}^{rs} } \, dz \Bigr )^{1/r}\ \ \leq\ C\ \norm
{V}_{L^{r'}}\ \norm {w}^2_{\domain{\Om}}\ \leq\ \frac12\  \norm
{w}^2_{\domain {\Om}},
\end{equation*}
hence $\norm {w}_{\domain{\Om}}=0$ and thus $w=0$, i.e., $u\leq \bar
u$.  The proof is complete.

\end{proof}
\medskip
\noindent In the proof of the previous Theorem we used  the
following comparison/uniqueness principle for the p-laplacian on an
unbounded domain.

\begin{lemma}\label{L:uniqueness p-lapl}
Let $V\in L^{r'}$ and $u,\ \bar u\in \dom$, $u\geq 0$,  be weak
solutions, respectively, of
\begin{align}\label{e:u bar eqn bis}
-\triangle_p \ u \ \leq \  V \frac {{u}^{p-1}}{\abs{x}^s}\ +\ g\\
\notag
\\ \label{eq for bar u}
-\triangle_p \ \bar u \ = \  V \frac {\abs{\bar u}^{p-2}\bar
u}{\abs{x}^s} \ +\ g.
\end{align}

a) Suppose $\bar u\geq 0$ and  $ g \geq 0$. If $g\not=0$ then $u\
\leq \ \bar u$. Otherwise, $u\ =\ c\, \bar u$ on the set $\{ u \geq
\bar u \}$.

b) More generally, suppose $g \,\bar u\geq 0$ and $\bar u$ is a
super solution of \eqref{eq for bar u}, i.e., $$ -\triangle_p \ \bar
u \ \geq \  V \frac {\abs{\bar u}^{p-2}\bar u}{\abs{x}^s} \ +\ g.$$
If $g \,\bar u\not\equiv 0$ then $u\ \leq \ \abs{\bar u}$.
Otherwise, $u\ =\ c\, \abs{\bar u}$ on the set $\{ u \geq \abs {\bar
u} \}$.

\end{lemma}

\medskip
\begin{proof}
a) For ease of reading let us consider first the case when equality
holds in the inequality satisfied for $u$, i.e.,  suppose $u,\ \bar
u\in \dom$ are two non-negative weak solutions. We are going to show
that if $g\not=0$ then $u\ =\ \bar u$. Otherwise, one of the
solutions is a constant multiple of the other. Working as in Lemma
3.1 in \cite{Lin}, see also \cite{DS} and \cite{A},  we define

$$
\ueps \ =\ u\ +\ \epsilon \quad \text { and }\quad \bueps \ =\  \bar
u \ +  \ \epsilon
$$
and the two test function
$$
\eta \ =\ \frac {\ueps^p - \bueps^p}{\ueps^{p}}\ =\ \ueps\ -\ \bigl
(\frac {\bueps}{\ueps} \bigr )^{p}\, \ueps
 \quad \text {
and}\quad \bar \eta \ =\  \frac {\bueps^p - \ueps^p}{\bueps ^{p}}\
=\ =\ \bueps\ -\ \bigl (\frac {\ueps}{\bueps} \bigr )^{p-1}\,
\bueps.
$$

\noindent Multiplying the equation for $u$ by $\eta$ and the
equation for $\bar u$ by $\bar \eta$, and then adding the two
equations we have

\begin{multline}\label{e:eta and eta bar}
\int_{\Om} \abs {\nabla u}^{p-2}\ \nabla u\cdot \nabla \eta\ +\ \abs
{\nabla \bar u}^{p-2}\ \nabla \bar u\cdot \nabla \bar\eta\, dz\\
 =\
\int_{\Om}  \frac {V}{\abs{x}^s} \ \bigl ( u^{p-1}\eta \ +\ {\bar
u}^{p-1} \bar \eta \bigr )\, dz \ +\ \int_{\Om} g \ ( \eta \ +\ \bar
\eta )\, dz.\hskip.4truein
\end{multline}

\noindent Using $$\nabla \eta= \Bigl [\, 1\ +\ (p-1) \Bigl ( \frac
{\bueps}{\ueps}\Bigr)^p \, \Bigr ]\ \nabla \ueps \ -\  p\ \Bigl (
\frac {\bueps}{\ueps}\Bigr)^{p-1} \ \nabla \bueps $$ we find

\begin{multline*}
\abs {\nabla u}^{p-2}\ \nabla u\cdot \nabla\eta\ =\ \Bigl [\, 1\ +\
(p-1) \Bigl ( \frac {\bueps}{\ueps}\Bigr)^p \, \Bigr ]\
\abs{\nabla\ueps}^p\ -\ p\ \Bigl ( \frac
{\bueps}{\ueps}\Bigr)^{p-1}\ \abs{\nabla\ueps}^{p-2}\
\nabla\ueps\cdot\nabla\bueps\\
\\
=\ \Bigl ( \ueps^p \ +\ (p-1)\ \bueps^p \Bigr )\ \abs{\nabla \ln
\ueps}^p\ -\ p\ \bueps^p\ \abs {\nabla\ln\ueps}^{p-2}\ \langle
\nabla\ln\ueps,\, \nabla \ln\bueps \rangle\\
\\
=\ \ueps^p \ \abs{\nabla \ln \ueps}^p\ +\ \bueps^p\ \Bigl [ -
\abs{\nabla \ln \ueps}^p\ -\ p \ \abs{\nabla \ln \ueps}^{p-2} \
\langle \nabla\ln\ueps,\, \nabla \ln\bueps - \nabla
\ln\ueps \rangle\Bigr ]\\
\\
\geq\ \ueps^p \ \abs{\nabla \ln \ueps}^p\ +\ \bueps^p\ \Bigl [\ C_p\
\abs {\nabla\ln \ueps \ -\ \nabla\ln\bueps }^p - \abs{\nabla \ln
\bueps}^p \ \Bigr ]\\
\\
=\ \ueps^p \ \abs{\nabla \ln \ueps}^p\ -\ \bueps^p \ \abs{\nabla \ln
\bueps}^p\ +\ C_p\ \bueps^p\ \abs {\nabla\ln \ueps \ -\
\nabla\ln\bueps }^p,\hskip.3truein
\end{multline*}
\noindent using the inequality, cf. Lemma 4.2 of \cite{Lin}, $$
\abs{a}^p\
> \ \abs{b}^p\ +\ p\ \abs {b}^{p-2}\ \langle b,\, a\ -\ b\rangle+ C_p\ \abs {a\ -\
b }^p\ , \quad a, \ b \in \Rn, \quad a\not= b.$$

\noindent The inequality we just proved shows the left-hand side of
\eqref{e:eta and eta bar} can be estimated as follows
\begin{multline*}
\int_{\Om} \abs {\nabla u}^{p-2}\ \nabla u\cdot \nabla \eta\ +\ \abs
{\nabla \bar u}^{p-2}\ \nabla \bar u\cdot \nabla \bar\eta\, dz\\
\geq \ C_p\ \int_{\Om} (\ueps^p \ + \ \bueps ^p)\ \abs {\nabla\ln
\ueps \ -\ \nabla\ln\bueps }^p\, dz.\hskip1truein
\end{multline*}

\noindent Therefore we have
\begin{multline*}
C_p \int (\ueps ^p + {\bueps} ^p)\ \abs {\nabla \ln \ueps
-\nabla\ln\bueps}^p \, dz \ \leq\ \int \frac {V}{\abs{x}^s} \Bigl [
  \Bigl ( \frac {u}{\ueps} \Bigr )^{p-1} \
   -\ \Bigl ( \frac {\bar u}{\bueps} \Bigr )^{p-1} \Bigr ]\ (\ueps^p - \bueps^p)\, dz
  \\
  \\
   +\ \int g\ \bueps\ \Bigl [ \frac {\ueps}{\bueps} + 1 - \Bigl ( \frac {\ueps}{\bueps}\Bigr)^p\ -
   \Bigl ( \frac {\ueps}{\bueps}\Bigr)^{1-p}\     \Bigr ]\,
   dz. \hskip.8truein
\end{multline*}
The first integral on the right goes to 0 when $\epsilon\rightarrow
0$ by the Lebesgue convergence theorem since $V\in L^{r'}$. A small
calculation shows that the function
\[
f(t)\ =\ 1\ +\ t \ - \ t^p \ -\  t^{1-p}
\]
is negative for $t>0$ with  $f=0$ iff $t=1$. In fact $f(1)=0$,
$f''(t)<0$ and $f'(1)=0$. The proof of the case when we have
equality in both places is complete.

Now, consider the more general case. Notice that $\eta \geq 0$ iff
$\bar \eta \leq 0$. Therefore we can multiply the inequality
satisfied by $u$ with the test function  $\eta^+$ and then add the
equation satisfied by $\bar u$ after we multiply it with $\bar\eta
^{-}\ =\ \min \{\bar\eta, 0 \}$. Then we work as before, but this
time the conclusions will hold only on the set $\{ u\geq \bar u\}$.
The proof of part (a) is complete.

b)  Define $v=\bar u^{+}=\max \{\bar u, 0 \}$, hence $0\leq v\in
\dom$, $\veps = v + \epsilon $ and
consider 

\begin{equation}
\bar \eta \ =\  \frac {\veps^p - \ueps^p}{\veps ^{p}}\ =\ =\ \veps\
-\ \bigl (\frac {\ueps}{\veps} \bigr )^{p}\, \veps.
\end{equation}

Therefore we can multiply the inequality satisfied by $u$ with the
test function  $\eta^+$ and then add the inequality satisfied by
$\bar u$ after we multiply it with $\bar\eta ^{\, -}\leq 0$, which
will bring us to
\begin{multline*}
C_p \int (\ueps ^p + {\veps} ^p)\ \abs {\nabla \ln \ueps
-\nabla\ln\veps}^p \, dz \ \leq\ \int \frac {V}{\abs{x}^s} \Bigl [
  \Bigl ( \frac {u}{\ueps} \Bigr )^{p-1} \
   -\ \Bigl ( \frac {v}{\veps} \Bigr )^{p-1} \Bigr ]\ (\ueps^p - \veps^p)\, dz\\
  \\
   +\ \int g\ \veps\ \, f({\ueps}/{\veps})
   dz. \hskip.8truein
\end{multline*}
Since $g\, \bar u \geq 0$, letting $\epsilon\rightarrow 0$ we see
that with $v=\bar u^+$ we have $u\leq v$ if $g\, \bar u >0$
somewhere, and $u$ and $v$ are proportional on the set $\{ u\geq
v\}$ otherwise, i.e., in the case $g\, \bar u\equiv 0$. Finally, let
us observe that $-\bar u$ is a solution of
$$ -\triangle_p \ (- \bar u) \ \geq \  V \frac {\abs{\bar u}^{p-2}(-\bar
u)}{\abs{x}^s} \ -\ g. $$ Taking $v=\max\{-\bar u, 0 \}= \ - \bar
u^{\, -}$ and observing that $-g \veps\ =\ g \bar u^{\, -}\ -\
\epsilon \, g$ the argument above shows that  $u\leq v$ if $g\, \bar
u >0$ somewhere, and $u$ and $v$ are proportional on the set $\{
u\geq v\}$ otherwise, i.e., in the case $g\, \bar u\equiv 0$. The
conclusion is that if $g\, \bar u \not\equiv 0$ then
\[
u\ \leq \ \abs{\bar u},
\]
while if $g\, \bar u \equiv 0$ then $u=c\bar u$ on the set $u\geq
\abs {\bar u}$.
\end{proof}

\medskip

So far we have concerned ourselves with the global properties of
solutions. Having done this we can obtain the asymptotic behavior of
solutions at infinity. The first result concern the local behavior
on a ball away from some finite point, see \cite{E} and \cite{Z} for
related results.

\medskip
\begin{thrm}\label{T:Asymptotic1}
Suppose $s,\ p, \ k,$ and $n$ satisfy the conditions of Theorem
\ref{boundedness}. Let $u$ be a nonnegative solution to the
inequality \eqref{e:maineq}, with $V\in L^{r'}$. We assume that $u$
has been extended with zero outside $\Om$. Suppose that $q_o\geq p$
is an exponent such that $u\in L^{q_o}(\Om)$. There exist constants
$C=C(\bG,p,||u||_{\mathcal{D}^{1,p}(\Om)}, \norm{u}_{q_o})>0$ and
$0<R_o=R_o(\norm{V}_{L^{r'}})$, such that, for every $z\in \bG$ and
$R=\abs{z}/2\geq R_o$ we have
\begin{equation}\label{locest}
\underset{B(z,R/2)}{max\ u} \leq C \left(\frac{1}{B(z,R)}
\int_{B(z,R)} u^{q_o} dx \right)^{\frac{1}{q_o}} .
\end{equation}
Furthermore,  $u$ has the following decay at infinity
\begin{equation}\label{e:asymp_decay}
u(z)\ \leq \ \frac {C}{\abs {z}^{n/{q_o}}}.
\end{equation}
\end{thrm}

\medskip

\begin{proof}[\textbf{Proof}]

Given a function $\alpha \in C^\infty_o(\bG)$, $\alpha \geq 0$ we
consider the function $\alpha^p F(u) \in \dom$, see \eqref{Fdef} and
\eqref{e:definition of G}. We recall \eqref{FGineq}, and observe in
addition
\begin{equation}\label{FGineq2}
F(u)\ \leq \ |G'(u)|^p.
\end{equation}
We use the fact $q/p>1$ in order to see that $G$ is a piece-wise
smooth and globally Lipschitz function. Using $\alpha^p F(u)$ as a
test function in the weak formulation \eqref{weakVsol} we have
\[
\int_{\Omega} \lvert {\nabla}u \rvert \,^{p-2} <{\nabla}u, {\nabla}
(\alpha^p F(u))> \,dz\ \leq \ \int_{\Omega} V\frac {\
|u|^{p-2}}{|x|^s}u\,\alpha^p F(u)\,dz
\]

Let us consider the left-hand side (LHS) of the above inequality,
which is easily seen to equal
\begin{equation}\label{e:LHS}
\text{LHS}\ =\ \int \alpha^p \abs{\nabla G}^p\, dz \ +\ p\int
\abs{\nabla u}^{p-2}\alpha^{p-1} F(u)\nabla u \cdot \nabla \alpha\,
dz.
\end{equation}
For any $\epsilon>0$ we have $ab\leq \epsilon \frac
{a^p}{p}+\epsilon^{-p'/p}\frac {b^{p'}}{p'}$, and hence
\begin{multline}\label{e:holder}
\abs{\nabla u}^{p-2}\alpha^{p-1} F(u)\nabla u \cdot \nabla \alpha\
\leq\ \frac {\epsilon }{p}\alpha^p\abs{\nabla u}^p u^{-1} F(u)\ +\
\frac {\epsilon^{-p'/p}}{p'}
\abs{\nabla\alpha}^p F(u)u^{p/{p'}}\\
\leq\ \frac {\epsilon }{p} \alpha^p \abs{\nabla G}^p \ +\
C\epsilon^{-1} q^{p-1} \abs {\nabla \alpha}^p G^p,
\end{multline}
after using  \eqref{FGineq2} and \eqref{FGineq} in the last
inequality.  Inserting \eqref{e:holder} in \eqref{e:LHS} we find
\[
\text{LHS}\ \geq\ (1-\epsilon)\int \alpha^p \abs{\nabla G}^p\, dz \
-\ C\epsilon^{-p'/p}q^{p-1}\int \abs {\nabla \alpha}^p G^p\, dz.
\]
We shall use the above inequality with a fixed sufficiently small
$\epsilon$ so that $1- {\epsilon }>0$. With the help of $(a+b)^p\leq
2^p(a^p+b^p)$, and using once more \eqref{FGineq} and the paragraph
after it, we come to
\begin{equation}
\int \abs {\nabla(\alpha G)}^p\, dz\ \leq\ C\int_{\Omega} V\frac
{\alpha^pG^p}{|x|^s}\,dz \ +\ C q^{p-1}\int \abs {\nabla \alpha}^p
G^p\, dz,
\end{equation}
where $C$ is a constant independent of $q$. Therefore, using the
H\"older and  the Hardy-Sobolev inequalities ( $pr=p^*(rs)$! ), we
have
\begin{equation}
\norm {\nabla(\alpha G)}^p_{L^p} \ \leq\ C \norm
{V}_{L^{r'}(\text{supp}\, \alpha)}\,\norm {\nabla(\alpha G)}^p_{L^p}
\ +\  C q^{p-1}\int \abs {\nabla \alpha}^p G^p\, dz
\end{equation}

Since $V\in L^{r'}$ it follows that if $0\not=z\in\Rn$ and
$R=\abs{z}/2\geq R_o$ then
\begin{equation}
\int_{B_R(z)} V^{r'} \, dz\ \rightarrow 0\quad\quad \text{as}
\quad\quad R_o\rightarrow\infty.
\end{equation}
Therefore, $C \norm {V}_{L^{r'}(\text{supp}\, \alpha)}\leq 1/2$ for
all $\alpha$ with

\begin{equation}\label{suppalpha}
 \alpha\in C^\infty_o\bigl ( B_R(z)\bigr ) \quad \text{ with }\quad R=\abs{z}/2, \quad\quad \abs{z}\geq R_o,
\end{equation}
 where $R_o$ depends on $V$ and $C$, and shall be fixed for the rest of the
proof. Using the Sobolev inequality, we have shown that for any such
$\alpha$ we have
\begin{equation}
\norm {\alpha G}_{L^{p^*}}\ \leq\ \norm {\nabla(\alpha G)}_{L^p} \
\leq\ C q^{(p-1)/p}\norm {(\nabla \alpha) G}_{L^p}
\end{equation}
Therefore, if $u\in L^q$ we can apply Fatou's theorem when $l
\rightarrow\infty$ to get
\begin{equation}
\Bigl ( \int \alpha^{p^*}u^{\delta q} \Bigr )^{1/{p^*}}\ \leq\ C\
q^{(p-1)/p} \Bigl ( \int \abs {\nabla \alpha}^p u^q  \Bigr )^{1/p},
\end{equation}
where $\delta=p^*/p>1$.  In particular, for any $0<\rho<r<R$ and
$\alpha\in C^{\infty}_o\bigl ( B(z,r)\bigr )$ with $\alpha\equiv 1$
on $ B(z,\rho)$ and $\abs {\nabla \alpha}\leq \frac {2}{r-\rho}$, we
have

\begin{equation}\label{e:moser}
\Bigl (\  \frac {1}{\abs{B(z,\rho)}}\int_{B(z,\rho)}u^{\delta q} \,
dz\ \Bigr)^{\frac {1}{\delta q}}\ \leq\  \frac
{C^{p/q}q^{(p-1)/q}}{(r-\rho)^{p/q}}\ \Bigl (\  \frac
{1}{\abs{B(z,r)}} \int_{B(z,r)}u^{ q} \, dz\  \Bigr )^{\frac {1}{
q}}.
\end{equation}
We can define the sequences $q_j=q_o\delta^j$ and $r_j=\frac
{R}{2}\bigl ( 1+\frac {1}{2^j}\bigr)$ for $j=0,1,\dots$ with which
Moser's iteration procedure gives inequality \eqref{locest}. Let us
observe that $\sum_{j=0}^{\infty} \frac {1}{q_j}<\infty$ and
$\sum_{j=0}^{\infty} \frac {\ln q_j}{q_j}<\infty$ thanks to
$\delta>1$, cf. \eqref{e:finitesums}. The decay property follows
immediately taking into account that the volume of $B_R(z)$ is
proportional to $R^n$, i.e., $\abs{z}^n$.
\end{proof}

Combining Theorems \ref{boundedness}, \ref{T:belowp*} and
\eqref{e:asymp_decay}, we can assert the following decay of
solutions to \eqref{e:maineq}.

\medskip

\begin{thrm}\label{T:asymptotic_beh p-lapl}
Let $\Om$ be an open subset of $\Rn$, which is not necessarily
bounded, $1< p <n$, $0 \leq s \leq p$, $s< k$ and $s(n-k)\ <\
k(n-p)$. Suppose $R\in L^{r'}\cap L^{t_o}$, for some $t_o >r'$, and
in the case $p>2$ assume that $R$ and $V_o$ are non-negative, $R\geq
0$, $V_o\geq 0$.

\noindent If   $u$ is a non-negative solution of
\eqref{e:main_ineq+Vo} then there exists a
$C=C(\bG,p,||u||_{\dom})>0$, such that,  $u$ has the following decay
at infinity
\begin{equation}
u(z)\ \leq \ \frac {C}{\abs {z}^{q}}\, \norm{u}_{D^{1,p}(\Omega)},
\end{equation}
for any $q<\frac {n-p}{p-1}$.
\end{thrm}

\begin{rmrk}
Let us observe that the fundamental solution of the p-laplacian on
$\Rn$ equals $C\abs {z}^{-\frac {n-p}{p-1}}$, where $C$ is a
constant.
\end{rmrk}

\vskip.6truein
\section{Asymptotics for the scalar curvature equation}
\vskip.2truein

In this section we restrict our considerations to the case $p=2$ and
furthermore we require $V=R(z)u^{2^*-2}$ with $R\in L^\infty$, i.e.,
we consider a non-negative weak solution $u$ of
\begin{equation}\label{e:Yamabetype}
-\triangle \, u\ \leq \ \frac {R(z)u^{2^*(s)-1}}{|x|^s} \quad \quad
\text{in} \quad \Omega.
\end{equation}

\noindent From Theorem \ref{T:asymptotic_beh p-lapl} we know that
for any $0<\theta<1$ there exists a constant $C_\theta>0$, such
that,
\begin{equation}\label{e:asymp_decay_laplacian}
u(z)\ \leq \ \frac {C_\theta}{1+\abs {z}^{\theta (n-2)}}\,
\norm{u}_{D^{1,p}(\Omega)}.
\end{equation}
The next result shows that the decay is at least as that of the
fundamental solution. Furthermore, if $u$ is a non-negative
solution, rather than a subsolution then $u$ has the same decay as
the fundamental solution.

\medskip
\begin{thrm}\label{T:asymptotic_beh lapl}
Let $\Om$ be an open subset (not necessarily bounded) of $\Rn$,
$n>2$, $0 \leq s \leq 2$, $s< k$ and  $R\in L^\infty$.

\noindent a) If $u$ is a non-negative solution of
\eqref{e:Yamabetype} then there exists a constant $C>0$, such that,
\begin{equation}\label{e:subsolutionest}
 0\ \leq\ u(z)\ \leq \ \frac
{C}{1+\abs {z}^{n-2}}, \quad z\in\Om.
\end{equation}
\noindent b) If $u$ is a non-negative non-trivial solution of
\begin{equation}\label{e:Yamabetype_solution}
-\triangle \, u\ = \ \frac {R(z)u^{2^*(s)-1}}{|x|^s} \quad \quad
\text{in} \quad \Omega,
\end{equation}
with $R\geq 0$ then
\begin{equation}\label{e:solutionest}
  \frac
{C^{-1}}{1+\abs {z}^{n-2}}\ \leq\ u(z)\ \leq \ \frac {C}{1+\abs
{z}^{n-2}}, \quad z\in\Om.
\end{equation}
\end{thrm}
\medskip

\begin{proof}
We shall prove first the estimate for $u$ from above.  Let us extend
$u$ as a function on $\Rn$ be setting it equal to zero outside of
$\Om$.  Let us define $f=\frac {|R(z)|u^{2^*-1}}{|x|^s}$. Thanks to
Theorem \ref{T:asymptotic_beh p-lapl} we have for any $0<\theta<1$
\begin{equation}\label{e:asympt_of_f}
f(z)\ \leq \ C_\theta \frac {|R(z)|}{1+\abs{z}^{\theta
(n+2-s)}}\frac {1}{\abs {x}^s}.
\end{equation}
Consider the function $v=\Gamma*f$, where $\Gamma$ is the positive
fundamental solution of the laplacian $\Gamma (z)=\frac
{1}{n(n-2)\omega_n}|z|^{2-n}$ ( recall $n>2$ ), where $\omega_n$ is
the volume of the unit ball in $\Rn$, so that $-\lap \Gamma =
\delta$. For points $z,\zeta$ in $\Rn$ we shall write $z=(x,y)\in
R^k\times\Rnmk$ and $\zeta=(\xi, \eta)\in R^k\times\Rnmk$. With this
notation and using \eqref{e:asympt_of_f} we have ( $C_n=\frac
{1}{n(n-2)\omega_n}$ )
\begin{align*}
v(z)\ & =\ C_n\int_{\Rn} \frac {f(\zeta)}{\abs{z-\zeta}^{n-2}}\, d\zeta\\
&\\
& =\  C_n \int_{|z-\zeta|\leq \frac {|z|}{2}} \frac {f(\zeta)}{\abs{z-\zeta}^{n-2}}\, d\zeta\  \\
\\ &\hskip.8truein +\ C_n \int_{|z-\zeta|>\frac {|z|}{2}} \frac {f(\zeta)}{\abs{z-\zeta}^{n-2}}\, d\zeta  \\
&\\
& \hskip.5truein \leq\ \frac {C_\theta\norm
{R}_\infty}{1+|z|^{\theta(n+2-s)} \ }\ \int_{|z-\zeta|\leq \frac
{|z|}{2}}\
\frac {1}{\abs{z-\zeta}^{n-2}}\frac {1}{|\xi|^s}\, d\zeta\\
&\\
 & \hskip.8truein +\ \frac {C_\theta \norm {R}_\infty}{1+|z|^{n-2} }\ \int_{\Rk} \int_{\mathbb{R}^{n-k}}  \frac {1}{(1+|\zeta|^
{\theta (n+2-s)} )}\frac {1}{ |\xi|^s}\, d\eta\, d\xi\\
&\\
& \overset {def}{\equiv}\ \frac {C_\theta\norm
{R}_\infty}{1+|z|^{\theta(n+2-s)} \ }\ I_1\ +\ \frac {C_\theta \norm
{R}_\infty}{1+|z|^{n-2} }\ I_2\ \leq\ \frac {C_\theta \norm
{R}_\infty}{1+|z|^{n-2} },
\end{align*}
from Lemmas \ref{l:I1 estimate} and \ref{l:Betafn}  ( see also the
Remark following Lemma \ref{l:I1 estimate} for the case $k=n$).

Going back to the bound from above for $u$ we see that the
difference $u-v$ is a subharmonic function, which goes to zero at
infinity and is equal to zero on $\partial\Om$. From the weak
maximum principle we can conclude
\[
u-v \ \leq\ 0,
\]
from which \eqref{e:subsolutionest}.

If $u$ is a solution,  we can take a ball $B$ centered at the origin
and a constant $C$, such that, the subharmonic function $u-C\,
\Gamma$ goes to zero at infinity, equals to zero on $\partial\Om$,
and is positive on the (compact) boundary $\partial B$. From the
weak maximal principle $u-C\, \Gamma$ is positive everywhere.  The
proof of the Theorem is complete.
\end{proof}

\medskip

\begin{rmrk}\label{r:Kelvinasymp} The second part of of Theorem \ref{T:asymptotic_beh
lapl}, i.e., when we assume in addition $R\geq 0$, can be derived
also with the help of the Kelvin transform on $\Rn$, which can be
seen as follows.
\end{rmrk}

 We know that the Kelvin transform is an isometry between $\domain
{\Om}$ and $\domain {\Om ^*}$, cf. \cite{E} and \cite{GV1}. A
calculation using  $\triangle\, (\mathcal{K} u)\ =\
|z|^{-n-2}(\triangle u)\, (\frac {z}{|z|^2})$ shows that the Kelvin
transform $\mathcal{K} u$ of $u$ is a non-negative weak solution of
the inequality
\[
-\triangle \,( \mathcal{K} u)(z)\ \leq\ R(\frac {z}{|z|^2})\, \frac
{|x|^s }{|z|^{2s}}(\mathcal{K} u)^{2^*(s)-1}(z)\ \leq\ R(\frac
{z}{|z|^2})\,\frac {1}{|x|^{s}}(\mathcal{K} u)^{2^*(s)-1} (z) \quad
\quad \text{in} \quad \Omega^*.
\]
Notice that $R(\frac {z}{|z|^2})\in L^\infty\, (\Om^*)$ when $R\in
L^\infty\, (\Om)$.  Furthermore, if $\Om$ is a neighborhood of the
infinity then $0$ is a removable singularity of $\mathcal{K} u$,
i.e., the above inequality is satisfied on $\Om^*\cup\{0\}$ and
$\mathcal{K}u\in\domain {\Om^*\cup\{0\}}$, cf. \cite{E} and
\cite{GV1}. In other words, when $0\leq R\in L^\infty$ the function
$Ku$ satisfies the same inequality and conditions as $u$. Recalling
that $u(z)={|z|}^{2-n}\, (\mathcal{K} u) (\frac {z}{|z|^2})$ the
claim follows from the regularity we have shown in Section 2.

\medskip
The following lemma contains a useful fact concerning the value of a
certain integral, which was used in the above Theorem.

\medskip
\begin{lemma}\label{l:Betafn}
Let  $m>\frac {n-s}{2}$ and $0\leq s <k < n$. The following formula
holds
\begin{multline}\label{e:Betafn}
\underset {\Rk}{\int} \underset {\Rnmk}{\int} \ \frac
{1}{(1+|x|^2+|y|^2)^m}\ \frac {1}{|x|^s} \, dy\,dx\\
 =\  \frac
{\sigma_{n-k}}{2}\frac {\sigma_{k}}{2}\ B(\frac {n-k}{2},m-\frac
{n-k}{2} )\ B(\frac {k-s}{2},m-\frac {n-s}{2} ).\quad
\end{multline}
\end{lemma}
\medskip

\begin{proof}
With $a^2=1+|x|^2$ we have
\begin{multline}
\int_{\Rnmk}\ \frac {1}{(1+|x|^2+|y|^2)^m}\, dy\  =\ \frac
{a^{n-k}}{a^{2m}} \int_{\Rnmk}\ \frac {1}{(1+|y|^2)^m}\, dy\\
\\
 =\ \frac {\sigma_{n-k}}{a^{2m-(n-k)}}\int_0^\infty\ \frac
{r^{n-k-1}}{(1+r^2)^m}\, dr\ =\ \frac
{\sigma_{n-k}}{2a^{2m-(n-k)}}\int_0^\infty\ \frac { t^{\frac
{n-k}{2}-1}}{(1+t)^{\frac {n-k}{2}+ (m-\frac {n-k}{2})}}\, dt\\
\\
 =\ \frac {\sigma_{n-k}}{2a^{2m-(n-k)}}\ B(\frac {n-k}{2},m-\frac
{n-k}{2} ).\quad\quad\quad
\end{multline}
Therefore we have
\begin{multline}
\int_{\Rk}\ \int_{\Rnmk}\ \frac {1}{(1+|x|^2+|y|^2)^m}\ \frac
{1}{|x|^s} \, dy\,dx\ \ \\
=\ \frac {\sigma_{n-k}}{2}\ B(\frac {n-k}{2},m-\frac {n-k}{2} )\
\int_{\Rk} \frac {1}{(1+|x|^2)^{m-\frac {n-k}{2}}}\
\frac {1}{|x|^s} \, dx\\
\\
= \ \frac {\sigma_{n-k}}{2}\frac {\sigma_{k}}{2}\ B(\frac
{n-k}{2},m-\frac {n-k}{2} )\ \int_0^\infty\ \frac { t^{\frac
{k-s}{2}-1}}{(1+t)^{\frac {k-s}{2}+ (m-\frac {n-k}{2}-\frac
{k-s}{2})}}\, dt\\
\\
= \ \frac {\sigma_{n-k}}{2}\frac {\sigma_{k}}{2}\ B(\frac
{n-k}{2},m-\frac {n-k}{2} )\ B(\frac {k-s}{2},m-\frac {n-s}{2} ).
\end{multline}

\end{proof}

\begin{rmrk}\label{R:Beta_Fn_Integral}
For future reference, let us notice that the above proof amounts to
using twice, with the appropriate choice of the involved parameters,
the formula
\[
\int_{\Rk} \frac {1}{(1+|x|^2)^a}\ \frac {1}{|x|^s} \, dx\ =\ \frac
{\sigma_{k}}{2}\ B(\frac {k-s}{2}, a - \frac {k-s}{2}),
\]
which is valid for any $a>0$, $k>s$ and $a >\frac {k-s}{2}$.
\end{rmrk}

We end the section with one more technical lemma, which was used in
Theorem \ref{T:asymptotic_beh lapl}.

\medskip
\begin{lemma}\label{l:I1 estimate}
Let $k\geq 2$, $n\geq 3$ and $0\leq\ s<k\ \leq\ n$ and $s<2$. There
exists a constant $C>0$ such that
\[
I(z)\ =\ \int_{|z-\zeta|\leq \frac {|z|}{2}}\ \frac
{1}{\abs{z-\zeta}^{n-2}}\frac {1}{|\xi|^s}\, d\zeta\ \leq \ C\
 |z|^{2-s}.
\]
\end{lemma}
\medskip

\begin{proof}
To estimate the integral $I$ we observe that $I$ is homogeneous,
$I(\lambda z)=\lambda^{2-s} I_1(z)$ for $\lambda>0$. Therefore, if
$I$ is finite on $\abs{z}=1$ we can conclude that
\begin{equation}\label{e:I1}
I(z)\ \leq \ C\abs{z}^{2-s}.
\end{equation}
In order to see that $I$ is finite when  $|z|=1$ let us notice that
it depends only on $|x|$ and $|y|$ as the integral is invariant
under rotation in $\Rk$ or $\Rnmk$. A consequence of this fact is
that it is enough to show that on $|z|=1$ the integral $I$ is finite
at only three point, namely, $x=0$, $x=1$,  and $x=y$ ( we write
$x=1$ for the number one on the real axis considered as a point in
$\Rk$, etc.).

For the rest of the proof we assume $|z|=1$.

 The case of $x\not=0$
is easier, so we shall consider the last two points first. Without
any loss of generality we take $x=1$ and we split the integral in
two parts
\[
I(z)\  =\ \int_{\{|z-\zeta|\leq \frac {|z|}{2}\}\cap \{|x-\xi|\leq
\frac {|z|}{4}\} \} } \quad\quad + \ \int_{\{|z-\zeta|\leq \frac
{|z|}{2}\}\cap \{|x-\xi|\geq \frac {|z|}{4}\} \} }.
\]

\noindent On the domain of integration of the first integral we have
that $|\xi|$ is bounded away from zero and so the integral is
finite. In turn, on the domain of integration of the second integral
$|z-\zeta|$ is bounded away from zero and $|\xi|^{-s}$ is integrable
near the origin as $k>s$, hence this integral is finite again.

Let us consider now the case $x=0$. Introducing $\rho\xi_o=\xi$ and
$r\eta_o=y-\eta$  we put $I$ in the form
\[
I\ =\ \int_{|\xi_0|=1}\, \int_{|\eta_o|=1}\ \int_{r^2+\rho^2\leq
1/4}\ \frac {\rho^{k-1}\ r^{n-k-1}}{(\rho^2+r^2)^{\frac {n-2}{2}}}\
\frac {1}{\rho^s}\, dr\, d\rho\, d\eta_o\, d\xi_o.
\]
Letting $r=t\cos \phi$, $\rho=t\sin \phi$ we come to
\[
I\ =\ \sigma_{n-k}\sigma_k \int_0 ^{1/2} t^{1-s}\, dt\ \int_0^{2\pi}
\frac {(\sin \phi)^{k-1} (\cos\phi)^{n-k-1}}{|\sin\phi|^s}\, d\phi\
< \infty,
\]
iff $s<2$ and $k>s$.
\end{proof}

\vskip .6in

\section{\textbf{A non-linear equation in $\mathbb{R}^n$ related to the Yamabe equation on groups of Heisenberg type}}
\vskip .2in

Suppose $a$ and $b$ are two natural numbers,  $\lambda>0$, and for
$x,y\in \mathbb{R}^{+}=(0,+\infty)$, define the function
\begin{equation*}\label{def of phi}
\phi\ =\ \lambda^2\bigl [\ (x+\alpha)^2 + (y +\beta)^2 \ \bigr],
\end{equation*}
where $\alpha,\ \beta\in \mathbb{R}$.
\medskip

\begin{prop}
The function $\phi$ satisfies the following equation in the plane
\begin{equation}
\Delta \phi\ - \ \frac{a+b+ 2}{2}\ \frac{|\nabla \phi|^2}{\phi}\ +\
\frac{a}{x}\ \phi_x\ +\ \frac{b}{y}\ \phi_y \ =\ \frac
{2a\lambda^2\alpha}{x}\ +\ \frac {2b\lambda^2\beta}{y}, \quad
xy\not=0.
\end{equation}
\end{prop}

\medskip

\begin{proof}
Set $\xi=\lambda (x+\alpha)$, $\eta=\lambda(y +\beta)$ and define
$\tphi(\xi,\eta) = \phi(x,y)$. Then we have
\[
\dx {} = \lambda\dxi {} \quad \text{ and } \quad \dy {} = \lambda
\deta {}.
\]
Thus we have,

\begin{align}
\Sigma\ \overset {def}{=} \  & \Delta \phi\ - \ \frac{a+b + 2}{2}\
\frac{|\nabla \phi|^2}{\phi}\ +\ \frac{a}{x}\ \phi_x\ +\
\frac{b}{y}\ \phi_y \notag\\
&\\
             & \quad\quad =\ \lambda\triangle\tphi\ - \ \frac{a+b + 2}{2}\lambda^2\ \frac{|\nabla \tphi|^2}{\tphi}\
              +\ \frac{a\lambda}{x}\
              \tphi_\xi\ +\
\frac{b\lambda}{y}\ \tphi_\eta.\notag
\end{align}
Since $\tphi  = \xi ^2 + \eta ^2$ we have

\begin{align}
\Sigma\  = \ & 4\lambda^2\ -\ \frac {n+2}{2}\lambda^2 \frac{4(\xi^2
+\eta^2)}{\xi^2 +\eta^2}\ +\ \frac ax2\lambda\xi\
+\ \frac by 2\lambda\eta\notag\\
&\notag\\
          &\quad =\ -2n\lambda^2\ +\ \frac {2a\lambda}{x}(\lambda x+\lambda\alpha)\
          +\ \frac {2b\lambda}{y}(\lambda y + \lambda\beta)\\
&\notag\\
          &\quad\quad\quad = \ -\ 2n\lambda^2 \ + \ 2a\lambda^2\ +\  2b\lambda^2 \ +\ \frac {2a\lambda^2\alpha}{x}\
          +\ \frac {2b\lambda^2\beta}{y}.\notag
\end{align}
Hence, taking into account $a+b=n$,  we proved
\[
\Sigma\ =\  \frac {2a\lambda^2\alpha}{x}\ +\ \frac
{2b\lambda^2\beta}{y}.
\]
\end{proof}

\medskip

Noting that $\Delta \phi \ +\ \frac{a}{x}\ \phi_x\  +\ \frac{b}{y}\
\phi_y $ is the laplacian in
$\mathbb{R}^n\equiv\mathbb{R}^{a+1}\times\mathbb{R}^{b+1}$,
$n=a+b+2$, acting on functions with cylindrical symmetry, i.e.,
depending on $|\vx|$ and $|\vy|$ only, we are lead to the following
question.

\medskip
\begin{Question} Given two real numbers $p_o$ and $q_o$ find all positive solutions of the equation
\[
\Delta u \ - \ \frac{n}{2}\ \frac{|\nabla u|^2}{u}\ =\ \frac
{p_o}{\abs {\vx}}\ +\ \frac {q_o}{\abs {\vy}}, \quad (\vx, \vy)\in
\mathbb{R}^n\equiv\mathbb{R}^{a+1}\times\mathbb{R}^{b+1},
\]
which have  at most a quadratic growth condition at infinity, $u\
\leq\ C (|\vx|^2 + |\vy|^2)$.
\end{Question}
\medskip

As usual a simple transformation allows to remove the appearance of
the gradient in the above equation. For a function $F$ we have
$\Delta F(u)\ =\ F''(u)|\nabla u|^2\ +\ F'(u)\Delta u$ and thus
\begin{align*}
\Delta u^\tau\ & =\ \tau (\tau -1) u^{\tau - 2} |\nabla u|^2 \ +\ \tau u^{\tau -1} \Delta u  \\
               & \quad\quad  =\ \frac {\tau}{2} u^{\tau -2}( 2u\Delta u\ +\ 2 (\tau-1)|\nabla u|^2).
\end{align*}
Therefore we choose $\tau$ such that $2(\tau -1)=-n$, i.e., $\tau\
=\ \frac {2-n}{2}$ and then rewrite the equation for $u$ as
\begin{align*}
\Delta u^{\frac {2-n}{2}}\ & =\ \frac {2-n}{2} u^{\frac {2-n}{2} -2} (2u\Delta u \ -\ n |\nabla u|^2)\label{eq for u}\\
                            & =\ - \frac {(n-2)}{2}u^{\frac {2-n}{2} -1}\ \Bigl ( \, \frac {p_o}{\abs {\vx}}\ +\ \frac {q_o}{\abs {\vy}}\, \Bigr ).
\end{align*}
This is the equation which we will study. As a consequence of the
above calculations we can write a three parameter family of
explicit solutions.
\medskip

\begin{prop}
Let $\lambda > 0$. The function $v(\vx, \vy)$ defined in
$\mathbb{R}^n\equiv\mathbb{R}^{a+1}\times\mathbb{R}^{b+1}$ by the
formula
\begin{equation}\label{def of v}
v(\vx,\vy)\ =\ \lambda^{2-n}(\ (|\vx| + \alpha)^2\ +\ (|\vy| +
\beta)^2  \ )^{\frac {2-n}{2}},\quad\quad (\vx, \vy)\in
 \mathbb{R}^n\equiv\mathbb{R}^{a+1}\times\mathbb{R}^{b+1},
\end{equation}
satisfies the equation
\begin{equation}
\Delta v \ =\ -\, v^{\frac {n}{n-2}}\ \Bigl (\, \frac {p}{|\vx|}\ +\
\frac{q}{|\vy|}\, \Bigr ),
\end{equation}
where
\begin{equation}
 {p}\ =\ \alpha\ {(n-2)\lambda^2\ a}, \quad\quad q\ =\ \beta\ {(n-2)\lambda^2\ b}.
\end{equation}
\end{prop}

\medskip

Let us observe that the above equation is invariant under rotations
in the $\vx$ or $\vy$ variables. Also if $v$ is a solution then a
simple calculations shows that for any $t\not=0$ the function $v_t
(\vx,\vy) \ =\ t^{(n-2)/2}\, v(t \vx, t \vy)$ is also a solution.

Another observation is that the same principle works if we split
$\mathbb{R}^n$ in more than two subspaces. For example, if we take
three subspaces we can consider the equation

\[
\Delta v \ =\  v^{\frac {n}{n-2}}f(|\vx|, |\vy|, |\vz|), \quad \quad
f(\vx,\vy,\vz)\ =\  \frac {p}{\abs {\vx}}\ +\ \frac{q}{\abs {\vy}} \
+\ \frac{r}{\abs\vz}
\]
and ask the question of finding all positive solutions with the same
behavior at infinity as the fundamental solution. Clearly the
function

\[
 v\ =\ \lambda^{2-n}(\ (|\vx| + \alpha)^2\ +\ (|\vy| + \beta)^2 \ +\ (|\vz| +\gamma)^2 \ )^{\frac {2-n}{2}},
\]

\noindent with the obvious choice of $\alpha, \ \beta$ and $\gamma$
is a solution.

\vskip .6in
\section{\textbf{The best constant and extremals of the Hardy-Sobolev inequality }}

\vskip .2in

In this Section we give the proof of Theorem
\ref{T:bestconstantHardySobolev}. It was proven in \cite{SSW} that
there are extremals with cylindrical symmetry, i.e., functions
depending only on $|x|$ and $|y|$ for which the inequality becomes
equality. On the other hand, it was shown in \cite{MS} that all
extremals of inequality \eqref{L2HS} have cylindrical symmetry after
a suitable translation in the $y$ variable, see also \cite{CW} and
\cite{LW} for some related results.

\medskip

\begin{thrm}[\cite{MS}]\label{MS}
If $u\in \domain {\Rn}$ is a function for which equality holds in
\eqref{L2HS} then
\begin{enumerate}[i)]
\item for any $y\in\Rmminusk$ the function $u(., y)$ is a radially symmetric decreasing function in $\Rk$; \item
there exists a $y_o\in\Rmminusk$ such that for all $x\in\Rk$ the function $u(x,. + y_o)$ is a radially symmetric
decreasing function on $\Rmminusk$.
\end{enumerate}
\end{thrm}
\medskip

 We turn to the proof of Theorem \ref{T:bestconstantHardySobolev}, in which we find the extremals and the best constant in
\eqref{L2HS} in the case $\sigma p_\sigma=1$, i.e., $s=1$ in Theorem
\ref {BT}.

 \begin{proof}\emph{(of Theorem \ref{T:bestconstantHardySobolev})}
By Theorem 2.1 and Theorem 2.5 of \cite{BT} there is a constant $K$
for which \eqref{sharpHS} holds and this constant is achieved, i.e.,
the equality is achieved. A small argument shows that a non-negative
extremal $u$ of the naturally associated variational problem $\inf
\int_{\Rm} |\nabla u|^2\,dz$ subject to the constraint
\begin{equation}\label{constraint}
\int_{\Rmminusk}\int_{\Rk} \frac {|u|^{\frac {2(n-1)}{n-2}}}{|x|}\,
dxdy\ = \ 1
\end{equation}
satisfies the Euler-Lagrange equation
\begin{equation}
\triangle u\ = \ - \frac {\Lambda}{|x|}\, u^{\frac {n}{n-2}},
\quad\quad u\in D^{1,2}\, (\Rm),
\end{equation}
where $\Lambda\ =\ K^{\frac {2(n-1)}{n-2}}$. From Theorem
\ref{T:L-infinity} and standard elliptic regularity results we can
see that $v$ is a $C^\infty$ function on $|x|\not=0$. Furthermore,
$\nabla u\in L_{\text{loc}}^\infty(\Rn)$ and $u$ is $C^\infty$
smooth in the $y$ variables. In particular $u\in
C_{\text{loc}}^{0,\alpha}(\Rn)$ for any $0<\alpha<1$. In order to
see these claims let $v=u_{y_j}$ for some $j$.  Hardy's inequality
shows $\frac {v}{|x|}\in L^q_{\text{loc}}$ for any $1<q<k$. From
elliptic regularity $v\in W^{2,q}_{\text{loc}}(\Rn)$ for any $1<q<k$
and hence the Sobolev embedding gives $v\in W^{1,\, \delta
q}_{\text{loc}}(\Rn)$, where $\delta = \frac {n}{n-q}>1$. After
finitely many iterations we see that $v\in
W^{1,2}_{\text{loc}}(\Rn)$ from which we can invoke Remark
\ref{r:localboundedness} to conclude $v\in
L^\infty_{\text{loc}}(\Rn)$. The same argument can be done for the
higher order derivatives in the $y$ variables. For the $x$
derivatives we argue similarly. We consider ( $v=u_{x_i} $ )
\[
\triangle v\ =\ -\frac {V}{|x|}v\ -\ \frac {V_o}{|x|} .
\]
We note that $V, \ V_o\ \in L_{\text{loc}}^\infty(\Rn)$ and it is
not hard to see that the proof of Theorem \ref{T:L-infinity}, cf.
also Remarks \ref{R:L-infinity for sols} and Remark
\ref{r:localboundedness}, allows to conclude $v \in
L_{\text{loc}}^\infty(\Rn)$, and hence $\nabla_x u\in
L_{\text{loc}}^\infty(\Rn)$ as claimed.

From Theorem 1.1 of \cite{MS} we can assume that $u$ is an extremal
with a cylindrical function after a suitable translation in the $y$
variable.  Thus we can assume that $u$ has cylindrical symmetry.
Introducing $\rho=|x|$, $r=|y|$ we have that $u$ is a function of
$\rho$ and $r$. We define $U(\rho,r)=u$ by restricting $u$ to two
lines through the origin-one in $\Rk$, the other in $\Rnmk$.  From
the regularity of $u$ it follows that $U$ is a smooth function of
$r$ for any fixed $\rho$. For any fixed $r$ it is a smooth function
of $\rho$ when $\rho\not=0$, and Lipschitz for any $\rho$.
Furthermore, in the first quadrant $\rho>0, r>0$  of the
$\rho\,r$-plane it satisfies the equation
\begin{equation}
\Delta U \ =\ -\, \frac {\Lambda}{\rho}\ U^{\frac {n}{n-2}}\ .
\end{equation}
Using the equation and the smoothness of $U$ in $r$ it is not hard
to see that $U$ has bounded first and second order derivatives on $
(\,(0,1)\times(0,1)\,)$, cf. Lemma \ref{l:ODE}.

Let $\phi(\rho,r)\ =\ U^{-\frac {2}{n-2}}$. The calculations of
Section 4 show that $\phi$ satisfies the following equation in the
plane
\begin{equation}
\Delta \phi\ - \ \frac{n}{2}\ \frac{|\nabla \phi|^2}{\phi}\ +\
\frac{a}{\rho}\ \phi_\rho\ +\ \frac{b}{r}\ \phi_r\ -\ \frac
{2\Lambda}{n - 2}\frac {1}{\rho}\ =\ 0,
\end{equation}
where $a=k-1$, $b=n-k-1$. Let $\mu > 0$ and consider $\tilde\phi\ =\
\mu^{-1}\phi$. Clearly $\tilde\phi$ is a solution of
\[
\Delta \tilde\phi\ - \ \frac{n}{2}\ \frac{|\nabla
\tilde\phi|^2}{\tilde\phi}\ +\ \frac{a}{\rho}\ \tilde\phi_\rho\ +\
\frac{b}{r}\ \tilde\phi_r\ -\ \frac {2\Lambda}{\mu(n - 2)}\frac
{1}{\rho}\ =\ 0.
\]
Let us choose $\mu$ such that $\frac {2\Lambda}{\mu(n - 2)}\ =\
\frac {n-2}{2}$, i.e.,
\[
\mu\ =\ \frac {4\Lambda}{(n-2)^2}.
\]
With this choice of $\mu$ we see that $\tilde\phi$ satisfies
equation (4.11) in \cite{GV}. Moreover, a small argument using the
homogeneity of the Kelvin transform shows it satisfies the
asymptotic behavior (4.37) of \cite{GV}, except the inequality
 for the  derivatives hold only on $|x|\not=0$. We can apply (4.40) of \cite{GV} by noticing that the
 integrals on the $\rho$ and $r$ axis vanish as $U$ has bounded first and second order derivatives in the punctured
 neighborhood of any point from the closed first quadrant, a fact which we observed above. Hence (4.43) of \cite{GV} after
setting $|A|=\lambda$ gives

\[
\tilde\phi\ =\ \lambda^2 \bigl [  \ (r+ \frac {n-2}{4a\lambda^2} )^2
+  s^2\ \bigr ],
\]
Recalling that $\phi \ =\  \mu \,\tilde\phi$ and the value of $\mu$ we come to
\[
\phi\ =\ \lambda^2 \,\frac {4\Lambda}{(n-2)^2}\, \bigl [  \ (r+
\frac {n-2}{4a\lambda^2} )^2 +  s^2\ \bigr ].
\]
This shows that $v$ must equal

\begin{align*}
v\ & =\  \lambda^{-(n-2)} \,\Bigl ( \frac {4}{(n-2)^2}\Bigr
)^{-\frac {n-2}{2}} \, \Lambda^{-\frac {n-2}{2}}\,
\Bigl [  \ (|x|+ \frac {n-2}{4a\lambda^2} )^2 + |y|^2\ \Bigr ]^{-\frac {n-2}{2}}\\
   & =\ \lambda^{-(n-2)} \,\Bigl ( \frac {n-2}{2}\Bigr )^{n-2} \, {K^{-(n-1)}}\,
\Bigl [  \ (|x|+ \frac {n-2}{4a\lambda^2} )^2 + |y|^2\ \Bigr
]^{-\frac {n-2}{2}}.
\end{align*}

The value of $K$ is determined by \eqref{constraint} after fixing
$\lambda$ arbitrarily, say $\lambda=1$, since the value of the
integral in \eqref{constraint} is independent of $\lambda$. With
this goal in mind we set $p= \frac {n-2}{4a}$ and compute the
integral

\begin{align}\notag
1\ =\ \int_{\Rmminusk\Rk} & \frac {1}{|x|} \Bigl [ \  \Bigl (\frac
{n-2}{2} \Bigr )^{n-2} \frac {1}{K^{n-1}}  \frac {1}{ \bigl [\,
(|x|+ p)^2+ |y|^2 \, \bigr ]^{\frac {n-2}{2}}}\ \Bigr ]^{\frac
{2(n-1)}{n-2}}\ dxdy\notag\\
&\label{integral}\\
 & =\  \frac {1}{K^{\frac {2(n-1)^2}{n-2}}} \Bigl (\frac {n-2}{2} \Bigr )^{2(n-1)} \
\int_{\Rmminusk\Rk} \frac {1}{|x|} \frac {1}{  \bigl [\,  (|x|+
p)^2+ |y|^2 \, \bigr ]^{n-1} }\ dxdy \notag
\end{align}
\medskip
 \noindent Let $a=|x|+p$. Then we compute

\begin{align}\notag
\int_{\Rmminusk} \frac {1}{(a^2 + |y|^2)^{n-1}}\, dy\ &  =\ \frac
{1}{a^{n+k-2}} \  \int_{\Rmminusk} \frac
{1}{(1+|y|^2)^{n-1}}\, dy\\
& \\
    & =\ \frac {\sigma_{n-k}}{2a^{n+k-2}}\ B(\frac {n-k}{2}, \frac {n+k}{2} -1)\notag,
\end{align}
where $\sigma_{n-k}$ is the volume of the unit $n-k$ dimensional
sphere and $B(., .)$ is the beta function. On the other hand after a
simple computation we find
\begin{align}\notag
\int_{\Rk} \frac {1}{|x| (|x|+p)^{n+k-2} } \ dx\ &  =\  \frac
{\sigma_k}{p^{n+k+1}}\
\int_0^\infty \frac {r^{k-2}}{(r+1)^{n+k-2}}\, dr\\
& \\
   &\quad\quad  =\  \frac {\sigma_k}{p^{n+k+1}} B(k-1,n-1)\notag.
\end{align}
Plugging in \eqref{integral} come to
\begin{align}\notag
K^{\frac {2(n-1)^2}{n-2}}\ & =\  \Bigl (\frac {n-2}{2} \Bigr
)^{2(n-1)}\frac {\sigma_{n-k}}{2}\ B(\frac
{n-k}{2}, \frac {n+k}{2} -1) \frac {\sigma_k}{p^{n+k+1}}\ B(k-1,n+k-1)\label{K}\\
&\\  & =\ 2^{2k+3}
(n-2)^{n-k-3}(k-1)^{n+k+1}\sigma_{n-k}\sigma_kB(\frac {n-k}{2},
\frac {n+k}{2} -1)B(k-1,n-1)\notag.
\end{align}

The proof is complete taking into account the allowed translations in the $y$ variable.
 \end{proof}

\medskip

In the above proof we used the following simple ODE lemma, which can
be proved by integrating the equation.

\medskip

\begin{lemma}\label{l:ODE}
Suppose $f$ is a smooth function on $\mathbb{R}\setminus \{0\}$,
which is also locally Lipschitz on $\mathbb{R}$, i.e., on any
compact interval there is a constant $L$, such that,
$|f(t')-f(t'')|\leq\ L\,|t'-t''|$ for any two points $t', t''$ on
this interval.  If $f$ satisfies the equation
\[
f''+\frac {k}{t}f'\ =\ \frac {a}{t}\ +\  {b},\quad \quad t>0,
\]
where $k$ is a constant $k>1$ and $a$, $b$ are
$L^\infty_{\text{loc}}$ functions, then $f$ has bounded first and
second order derivatives near the origin.
\end{lemma}

\vskip .6in

\section{\textbf{Some applications }}

Let us consider the prescribed scalar curvature equation on $\Rn$
\begin{equation}\label{e:scalar curv}
\lap u \ = \ -R(z) \ u^{2^*-1},
\end{equation}
where $R$ is a bounded  function and $u$ is a non-negative function.
As usual we say that $u$ is of finite energy if $\norm {u}_\domG$ is
finite. Clearly if $g=u^{4/(n-2)}g_o$ is a metric conformal to the
Euclidean metric $g_o$ on $\Rn$ then the finite energy condition is
equivalent to $g$ having finite volume.  The results of Section 2
can be extended to the case of the scalar curvature equation of many
non-compact manifold with positive Yamabe invariant, which among
other things will be done in \cite{VZ}, but the following result,
which follows from the fast decay Theorem \ref{T:asymptotic_beh
lapl} (a) of $u$, i.e., at least as fast as the fundamental
solution, is indicative of what is to be expected, see also
\cite{Le}.

\begin{thrm}\label{T:finite energy sol for Yamabe}
Suppose $R\in L^\infty$  and $g_o$ is the Euclidean metric on $\Rn$.
Let $u$ be a positive solution to \eqref{e:scalar curv}. If $\Rn$
with the conformal metric $g=u^{4/(n-2)}g_o$ has finite volume then
$u$ has fast decay and the metric $g$ is incomplete.
\end{thrm}

\medskip
The second application concerns the original motivation of Badiale
and Tarantello \cite{BT} to consider the Hardy-Sobolev inequality.
The following equation has been proposed, cf. \cite{Ch}, \cite{B}
and \cite{R} for further details, as a model to study elliptic
galaxies.
\begin{equation}\label{e:galaxies}
-\triangle u\ =\ \phi(|x|) u^{q-1},\quad \ 0\ <\ u\in
D^{1,2}(\mathbb{R}^3)\quad (\ z=(x,y)\in \mathbb{R}^3 \ !).
\end{equation}

\noindent It is also required that $u$ is of finite mass, i.e.,
\[
\int \phi u^{q-1}\, dz\ <\ \infty.
\]
Using the results of this paper we can show the following Theorem.

\medskip

\begin{thrm}\label{T:galaxies}
Suppose $ (1+|x|)^\gamma\phi\in L^\infty (\mathbb{R}^3)$ for some
$0<\gamma<2$. If  $ 2^* (\gamma) <q<6$, then any solution of
\eqref{e:galaxies} is also of finite mass.
\end{thrm}
\medskip

\begin{proof}
Since the dimension of the ambient space is three we have $2^*(0)=6$
and  $2<2^*(\gamma)<6$. Given any $q$ satisfying $2^*(\gamma)<q< 6$
we can find an $s<\gamma$, such that, $q=2^*(s)$ and $u$ satisfies
the equation
\[
-\triangle u\ =\ \phi(|x|)u^{q-1}\ =\ V\frac {|u|^{2^*(s)-1}}{|x|^s}
\]
with $|V|=|x|^s|\phi|\leq (1+|x|)^s |\phi |\ \leq (1+|x|)^\gamma
|\phi|\in L^\infty (\mathbb{R})$. Theorem \ref{T:asymptotic_beh
lapl} (a) implies that $u$ decays at least as fast as the
fundamental solution of the laplacian in $\mathbb{R}^3$
\[
u(z)\ \leq\ \frac {C}{1+|z|}.
\]
Since $q-1>2^*(\gamma)-1>1$ it follows
\[
\int \phi u^{q-1}\, dz\ <\ \infty,
\]
which shows that every finite energy solution is also of  finite
mass.
\end{proof}

\end{document}